\documentclass[10pt,psamsfonts,leqno]{siamltex}

\usepackage{amsfonts,amssymb,amsmath,amsxtra,euscript}

\newtheorem{prop}{Proposition}[section]
\newtheorem{thm}[prop]{Theorem}
\newtheorem{cor}[prop]{Corollary}
\newtheorem{lem}[prop]{Lemma}
\newtheorem{defn}[prop]{Definition}
\newtheorem{example}[prop]{Example}

\newtheorem{rmks}[prop]{Remarks}
\newcommand{\bbC}{\mathbb{C}}

\newcommand{\bbR}{\mathbb{R}}
\newcommand{\bbT}{\mathbb{T}}
\newcommand{\bbZ}{\mathbb{Z}}

\newcommand{\bbL}{\mathbb{L}}
\newcommand{\bbG}{\mathbb{G}}

\newcommand{\calA}{\mathcal{A}}
\newcommand{\calB}{\mathcal{B}}
\newcommand{\calC}{\mathcal{C}}
\newcommand{\calD}{\mathcal{D}}

\newcommand{\calF}{\mathcal{F}}

\newcommand{\calH}{\mathcal{H}}

\newcommand{\calK}{\mathcal{K}}

\newcommand{\calS}{\mathcal{S}}

\newcommand{\calU}{\mathcal{U}}

\newcommand{\scrC}{\EuScript{C}}

\newcommand{\scrG}{\mathfrak{G}}
\newcommand{\scrF}{\mathfrak{F}}

\renewcommand{\mod}{\operatorname{mod}}

\newcommand{\Dom}{\operatorname{Dom}}
\newcommand{\esssup}{\operatornamewithlimits{ess\, sup}}

\renewcommand{\t}{\tau}

\newcommand{\vphi}{\varphi}

\newcommand{\disp}{\displaystyle}
\newcommand{\lb}{\label}

\newcommand{\CoR}{{C_0(\bbR,X)}}

\newcommand{\CoRp}{{C_0(\bbR_+,X)}}
\newcommand{\Coo}{{C_{00}(\bbR_+,X)}}

\newcommand{\LpR}{{L^p(\bbR,X)}}
\newcommand{\LpRp}{{L^p(\bbR_+,X)}}
\newcommand{\LpRX}{{L^p(\bbR_+,X)}}
\newcommand{\LpRU}{{L^p(\bbR_+,U)}}
\newcommand{\LpRY}{{L^p(\bbR_+,Y)}}

\newcommand{\Uone}{\mathcal{U}_1}

\newcommand{\U}{\mathcal{U}}

\newcommand{\rstab}{r_{stab}}
\newcommand{\rcstab}{rc_{stab}}

\newcommand{\Gammaper}{\Gamma_{per}}
\newcommand{\Gammaperxi}{\Gamma_{per,\xi}}
\newcommand{\Fper}{{\mathcal{F}_{per}}}

\renewcommand{\theenumi}{\roman{enumi}}

\title{Stability Radius and Internal versus
External\\ Stability in Banach Spaces:\\
an Evolution Semigroup Approach
\thanks{A portion of this article
appeared in preliminary form in the {\em Proceedings of the 36th IEEE
Conference on
Decision and Control}, 1997. }}

\author{Stephen Clark
\thanks{Department of Mathematics and Statistics, University of
Missouri--Rolla, Rolla MO 65409 ({\tt sclark@umr.edu}).}
\and Yuri Latushkin
\thanks{Mathematics Department, University of
Missouri--Columbia, Columbia, MO 65201
({\tt yuri@math.missouri.edu}).
Research of this author was supported by the National Science
Foundation and the Missouri Research Board.}
\and Stephen Montgomery-Smith
\thanks{Mathematics Department,
University of Missouri--Columbia, Columbia, MO 65201
({\tt stephen@math.missouri.edu}) Research of this author was
supported by the National Science Foundation.}
\and Timothy Randolph
\thanks{Department of Mathematics
and Statistics, University of Missouri--Rolla, Rolla MO 65409
({\tt randolph@umr.edu}). Research of this author was supported by
the Missouri Research Board; a portion of the work by this author
was carried out while visiting the Department of Mathematics,
University of Missouri--Columbia.}}

\begin{document}
\maketitle

\begin{abstract}
In this paper the theory of evolution semigroups is developed and
used to provide a framework  to study the stability of general linear
control
systems. These include  autonomous and nonautonomous systems modeled
with unbounded
state-space operators acting on Banach spaces.  This approach allows
 one to apply the classical theory of strongly continuous semigroups to
time-varying  systems.  In particular, the complex stability radius may
be expressed explicitly in terms of the generator of a
(evolution) semigroup.  Examples are given to show that classical
formulas  for the stability radius of an autonomous Hilbert-space system
fail in more general settings. Upper and lower bounds on the stability
radius are proven for  Banach-space systems.  In addition, it is
shown that the theory of evolution semigroups allows for a
straightforward
operator-theoretic analysis of internal stability as determined by
classical frequency-domain and input-output operators, even for
nonautonomous Banach-space systems. In particular, for the nonautonomous
setting,
internal stability is shown to be equivalent to input-output stability,
for stabilizable and detectable systems.  For the autonomous setting,
an explicit formula for the norm of input-output operator is given.
\end{abstract}



\pagestyle{myheadings}
\thispagestyle{empty}
\markboth{S.~CLARK, Y.~LATUSHKIN, S.~MONTGOMERY-SMITH, AND T.~RANDOLPH}
{STABILITY IN BANACH SPACES}

%\markboth{T.RANDOLPH, Y.LATUSHKIN ET.AL.}{STABILITY IN BANACH SPACES}
%\tableofcontents
%-------------------------- section 1 -----------------------
\section{Introduction}

Presented here is a study of stability of infinite-dimensional linear
control systems which is based on the relatively recent development of
the theory of evolution semigroups.  These semigroups have been used in
the study of exponential dichotomy of time-varying differential
equations and more general hyperbolic dynamical systems; see
\cite{ChiLat,LMS2,LMSR,MiRaSc,vanNbook,RS1,RolandDis} and the
bibliographies therein. The intent of this paper is to show how the
theory of evolution  semigroups can be used
to  provide a clarifying perspective, and prove new results, on the
uniform
exponential stability  for  general linear control systems,
$\dot x(t)=A(t)x(t)+B(t)u(t)$, \ $y(t)=C(t)x(t)$, \ $t\ge0$.
The operators
$A(t)$ are generally unbounded operators on a Banach space,
$X$, while the operators $B(t)$ and $C(t)$ may act on Banach spaces,
$U$ and $Y$, respectively.  In addressing
the general settings, difficulties arise both from the time-varying
aspect and from a loss of Hilbert-space properties.
This presentation, however, provides some relatively simple
operator-theoretic arguments for properties that extend classical
theorems of autonomous systems in finite dimensions.
  The topics covered here include characterizing internal stability of
the nominal system in terms of  appropriate input-state-output operators
and, subsequently, using these properties to obtain  new explicit
formulas
for bounds on the stability radius.
Nonautonomous systems are generally considered, but some results apply
only to autonomous ones, such as the upper bound for the stability
radius (Section~3.3),  the formula for the norm of the input-output
operator in Banach spaces (Section~3.4) and a characterization of
stability that is related to this formula (Section~4.2).

Although practical considerations usually dictate that $U$ and $Y$ are
Hilbert spaces (indeed, finite dimensional), the Banach-space
setting addressed here may be motivated by the
problem of determining optimal sensor (or actuator) location.
For this, it may
be natural to consider $U=X$ and $B=I_X$ (or
$Y=X$ and $C=I_X$) \cite{burns};
if the natural state space $X$ is a
Banach space then, as will be shown in this paper, Hilbert-space
characterizations of internal stability  or
its robustness do not apply.  We also show that even in
the case of Hilbert spaces $U$ and $Y$, known formulas for the
stability radius
involving the spaces $L^2(\bbR_+,U)$ and $L^2(\bbR_+,Y)$ do
not apply if the $L^2$ norm is replaced by, say, the $L^1$ norm---see
the  examples in Subsection~\ref{counterex} below.
In addition to the general setting of
nonautonomous systems on Banach spaces, autonomous and Hilbert-space
systems are considered.

For the autonomous case,  the primary observation we make
about general Banach-space settings versus the classical
$L^2$ and Hilbert-space setting
can be explained using the notion of $L^p$--Fourier multipliers.
For this, let $H(s)=C(A-is)^{-1}B$,
$s\in{\mathbb R}$, denote the transfer function and ${\mathbf F}$
denote the Fourier transform. The transfer function $H$ is said to be
an {\em $L^p$--Fourier multiplier} if the operator $u\mapsto {\mathbf
F}^{-1}H(\cdot){\mathbf F}u$ can be extended from the Schwartz class of
rapidly decaying $U$--valued functions to a bounded operator from
$L^p({\mathbb R};U)$ to $L^p({\mathbb R};Y)$;
see, e.g., \cite{Am2} for the definitions.
As shown in Theorem \ref{LnormThm} below,
the norm of this operator is equal to the norm
of the input-output operator.
If $U$ and $Y$ are Hilbert spaces and $p=2$, then $H$ is an $L^2$--Fourier
multiplier if and only if $\|H(\cdot)\|$ is bounded on ${\mathbb R}$;
see formula \eqref{p=2}. For Banach spaces and/or $p\neq 2$ this latter
condition is necessary, but {\em not sufficient} for $H$ to be an
$L^p$--Fourier multiplier.  As a result, our formula \eqref{Lnorm1} for
the norm of the input-output operator is more involved.

To motivate the methods,
 recall Lyapunov's stability theorem which says that if
$A$ is a bounded linear operator on $X$ and if the spectrum of $A$ is
contained in the open left half of the complex plane, then the solution
of the autonomous differential equation
$\dot{x}(t)=Ax(t)$ on $X$  is uniformly exponentially stable;
equivalently,  the spectrum $\sigma(e^{tA})$ is contained in the open
unit disk $\{z\in\bbC:|z|<1\}$, for $t>0$.
This is a consequence of the fact that when $A$ is a bounded operator,
then the spectral mapping theorem holds:
$
\sigma(e^{tA})\setminus\{0\}=e^{t\sigma(A)},\quad t\ne 0.
$
Difficulties with Lyapunov's Theorem arise when the operators, $A$, are
allowed to be unbounded.   In particular, it is well known that there
exist strongly continuous semigroups
$\{e^{tA}\}_{t\ge 0}$  that are {\em not} uniformly exponentially stable
even though $\mbox{Re}\lambda\le \omega<0$ for all
$\lambda\in\sigma(A)$; see, e.g., \cite{Nagel,vanNbook,Renardy}.
  For nonautonomous equations the situation is worse. Indeed,  even for
finite-dimensional $X$ it is possible for the spectra of $A(t)$ to be
the same for all $t>0$ {\em and} contained in the open left half-plane
yet the corresponding solutions to
$\dot{x}(t)=A(t)x(t)$ are not uniformly
exponentially stable (see \cite[Exm.7.1]{Hale} for a classical example).
In the development that follows we plan to show how these difficulties
can be overcome by the construction of an ``evolution semigroup."  This
is a family of operators defined on a  superspace  of functions from
$\bbR$ into $X$, such as
$\LpR$, $1\le p<\infty$, or
$\CoR$.

Section~1 sets up the notation and provides  background information.
Section~2 presents the basic properties of the evolution semigroups.
Included here is the property that the spectral mapping theorem always
holds for these semigroups when they are defined on $X$-valued functions
on the
{\em half-line}, such as $\LpRp$.  A consequence of this is a
characterization
of exponential stability for nonautonomous systems in terms of the
invertibility of the generator $\Gamma$ of the evolution semigroup.
This
operator, and its role in determining exponential stability, is the
basis for many of the subsequent developments. In particular, the
semigroup
$\{e^{tA}\}_{t\ge 0}$ {\em is} uniformly exponentially stable provided
$\text{Re }\lambda<0$ for all $\lambda\in\sigma(\Gamma)$.

Section~3 addresses the topic of the (complex) stability radius; that
is,
the size of the smallest disturbance, $\Delta(\cdot)$, under which the
perturbation, $\dot x(t)=(A(t)+\Delta(t))x(t)$, of an exponentially
stable system, $\dot x(t)=A(t)x(t)$, looses exponential  stability.
Results address structured and unstructured perturbations of autonomous
and nonautonomous systems in both Banach and Hilbert space settings.
Examples are given which highlight some important differences
between these settings.  Also included in this section is a
discussion about the transfer function for infinite-dimensional
{\em time-varying} systems.  This concept arises naturally in the
context of evolution semigroups.

In Section~4 the explicit relationship between internal and
external stability is studied for general linear systems.
This material  expands on the ideas begun in \cite{CDC}.
 A classical result for autonomous systems in
Hilbert space is the fact that  exponential stability of the nominal
system (internal stability) is, under the hypotheses of
stabilizability and detectability, equivalent to the boundedness of the
transfer function in the right half-plane (external stability). Such a
result does not apply to nonautonomous systems and a counterexample
shows
that this property fails to hold for Banach space systems.  Properties
from Section~2 provide a natural Banach-space extension of this result:
the role of transfer function is replaced by the input-output operator.
Moreover, for  autonomous systems we provide an explicit formula
relating the norm of this input-output operator to  that of the
transfer function.  Finally, we prove two theorems---one for
nonautonomous and one for autonomous systems---which characterize
internal stability in terms of the various input-state-output operators.

This Introduction concludes with a brief synopsis of the main results.
The characterization of uniform exponential stability
in terms of an evolution semigroup and its generator
is given in Theorem \ref{ExpStabThm+}, Theorem \ref{TVvanN}, and
Corollary \ref{vanNCor}.  Although these results are essentially known,
the proofs are approached in a new way.  In particular,
Theorem~\ref{TVvanN} identifies the operator $\bbG=-\Gamma^{-1}$  used
to determine
stability throughout the paper.
Theorem~ \ref{L0BoundThm} records the
 main observation that the input-output operator,
$\bbL=\calC\bbG\calB$, for a general nonautonomous system is related
to the inverse of the generator of the evolution semigroup.
A very short proof of the known fact that the stability radius
for such a system
is bounded from below by $\|{\mathbb L}\|^{-1}$ is also provided here.
The upper bound for the stability radius,
being given in terms of the transfer function, applies only
to autonomous systems and is proven in Subsection~3.3.
The upper bound, as identified here for Banach spaces,
seems to be new although
our proof is based on the idea of the
Hilbert-space result of \cite[Thm.~3.5]{HP94}.
In Subsection~3.3 we also
introduce the  pointwise stability radius and
dichotomy radius.  Estimates for the former are provided by
 Theorems \ref{LMS-Gearhart+BDeltaC} and \ref{PtwiseBounds} while the
latter is addressed in Lemma~\ref{dichrad}.
Examples \ref{RstabIneqEx} and \ref{StephenEx} show that,
for autonomous Banach space systems, both
inequalities for the upper and lower bounds on the stability radius
(see Theorem~\ref{autin}) can be strict.
In view of the
%% yuri-change
possibility of the
strict inequality
$\|\bbL\|>\sup_{s\in\bbR}\|C(A-is)^{-1}B\|$,
Theorem \ref{LnormThm} provides a new Banach space formula for
$\|\bbL\|$ in terms of $A$, $B$, and $C$.
In  Section~4 this expression
for $\|L\|$ is used to relate state-space versus frequency-domain
stability---concepts which are {\em not} equivalent for Banach-space
systems.
A special case of this expression gives a new formula
for the growth bound of a semigroup on a Banach space; see
Theorem~\ref{autonstab} and the subsequent paragraph.
Finally, Theorem~\ref{bigIOstab}  extends a classical
characterization of stability for stabilizable and detectable control
systems as it applies to nonautonomous Banach-space settings.

%------------------section 1--------------------

\section{Notation and Preliminaries}

Throughout the paper, $\mathcal{L}(X,Y)$ will denote the set of bounded
linear operators between complex Banach spaces $X$ and $Y$.
If $A$ is a linear operator on $X$, $\sigma(A)$ will denote the
spectrum of $A$, $\rho(A)$ the
resolvent set of $A$ relative to $\mathcal{L}(X)=\mathcal{L}(X,X)$, and
$\| A \| _\bullet=\| A \|_{\bullet,X}:=
 \inf \{\| Ax \|: x\in \Dom(A),  \|x \| =1\}$.
In particular, if $A$ is invertible in $\mathcal{L}(X)$,
$\| A\|_\bullet = 1/ \|A^{-1}\|_{\mathcal{L}(X)}$.
Also, let ${\mathbb C}_+=\{\lambda\in{\mathbb C}:
\text{Re }\lambda>0\}$.

If $A$ generates a strongly continuous (or $C_0$) semigroup
$\{e^{tA}\}_{t\ge 0}$
on a Banach space $X$ the following notation will be used:
$s(A)=\sup\{\text{Re}\lambda: \lambda\in\sigma(A)\}$ denotes the
spectral bound; $s_0(A)=
\inf\{\omega\in{\mathbb R}:\,\,\{\lambda:\text{Re
}\lambda>\omega\}\subset\rho(A)\quad\text{and}\quad\sup_{\text{Re
}\lambda>\omega}\|(A-\lambda)^{-1}\|<\infty\}$ is the abscissa of
uniform
boundedness of the resolvent; and
$\omega_0(e^{tA})=\inf\{\omega\in{\mathbb R}:\,\,\|e^{tA}\|\le
Me^{t\omega}\quad\text{for some}\quad M\ge 0\quad\text{and all}\quad
t\ge
0\}$ denotes the growth bound of the semigroup. In general,
 $s(A)\le s_0(A)\le \omega_0(e^{tA})$ (see, e.g., \cite{vanNbook})
 with {\em strict} inequalities possible; see
\cite{Nagel,vanNbook,Renardy}
for examples.  However, when $X$ is a Hilbert space, the following
spectral
mapping theorem of L.~Gearhart holds (see, e.g., \cite[p.~95]{Nagel} or
\cite{vanNbook,Prus}):

\begin{thm}\label{GearTh} If $A$  generates a strongly continuous
semigroup
$\{e^{tA}\}_{t\ge0}$ on a Hilbert space, then $s_0(A)=
\omega_0(e^{tA})$. Moreover, $1\in \rho(e^{2\pi A})$ if and only if
$i\bbZ \subset\rho (A)$ and
$\sup_{k\in \bbZ} \| (A-ik)^{-1}\|< \infty$.
\end{thm}

\noindent In particular, this result shows that on a Hilbert space $X$
the semigroup $\{e^{tA}\}_{t\ge 0}$ is uniformly exponentially stable if
and
only if $\sup_{\lambda\in{\mathbb C}_+}\|(A-\lambda)^{-1}\|<\infty$
\cite{Huang}.

Now consider operators $A(t)$, $t\ge0$, with domain $\Dom(A(t))$
in a Banach space $X$. If the abstract Cauchy problem
\begin{equation}\lb{nacp}
\dot x(t)=A(t)x(t),\quad x(\tau)\in\Dom(A(\tau)),\qquad t\ge \tau \ge 0,
\end{equation}
is well-posed in the sense that there exists
an evolution (solving) family of operators $\U=\{U(t,\tau)\}_{t\ge
\tau}$ on $X$ which gives a differentiable solution, then
$x(\cdot):t\mapsto U(t,\tau)x(\tau)$,
$t\ge \tau$ in $\bbR$, is differentiable, $x(t)$ is in $\Dom(A(t))$ for
$t\ge0$, and \eqref{nacp} holds.  The precise meaning of the term
evolution family used here is as follows.

\begin{defn}\lb{DEFevolfam}  A family of bounded operators
$\{U(t,\tau)\}_{t \geq \tau}$ on $X$
is called an {\em evolution family} if
\begin{enumerate}
\item
$U(t,\tau )=U(t,s)U(s,\tau )$ and $U(t,t)=I$ for all
$t \geq s\geq \tau$;
\item for each $x\in X$ the function $(t ,\tau)\mapsto
U(t,\tau)x$ is continuous for $t \geq \tau$.
\end{enumerate}
An evolution family $\{U(t ,\tau)\}_{t \geq \tau}$ is
called {\em exponentially bounded} if, in addition,
\begin{enumerate}
\item[(iii)]
there exist constants $M\geq 1$, $\omega\in\bbR$ such that
$$
\|U(t, \tau)\|\leq Me^{\omega (t -\tau)},\quad  t \geq \tau.
$$
\end{enumerate}
\end{defn}

\begin{rmks}
\begin{enumerate}
\item[(a)] An evolution family $\{U(t,\tau)\}_{t\ge\tau}$
is called {\em uniformly
exponentially stable}  if in part (iii), $\omega$ can be taken to be
strictly less than zero.
\item[(b)] Evolution families appear as solutions for abstract Cauchy
problems \eqref{nacp}.  Since the definition requires
 that $(t ,\tau)\mapsto U(t,\tau)$ is merely strongly continuous
the operators $A(t)$ in \eqref{nacp} can be unbounded.
\item[(c)] In the autonomous case where $A(t)\equiv A$ is the
infinitesimal generator of a strongly continuous semigroup
$\{e^{tA}\}_{t\ge 0}$ on $X$ then $U(t, \tau)=e^{(t-\tau )A}$, for
$t \geq \tau$, is a strongly continuous exponentially bounded
evolution family.
\item[(d)]  The existence of a {\em differentiable} solution to
\eqref{nacp} plays little role in this paper, so the starting point
will usually not be the equation \eqref{nacp}, but rather the
existence of an exponentially bounded evolution family.
\end{enumerate}
\end{rmks}

In the next section we will define the  evolution semigroup relevant to
our interests for the nonautonomous Cauchy problem \eqref{nacp} on the
half-line, $\bbR_+=[0,\infty)$.  For now,  we begin by considering the
autonomous equation $\dot x(t)=Ax(t)$, $t\in\bbR$, where $A$ is the
generator of a strongly continuous semigroup  $\{e^{tA}\}_{t\ge 0}$ on
$X$.  If $\mathcal{F}_\bbR$ is a space of $X$-valued functions,
$f:\bbR\to X$, define
\begin{equation}\lb{AUTevolsgR}
(E_\bbR^tf)(\t)=e^{tA}f(\t-t),\quad \mbox{for }f\in \mathcal{F}_\bbR.
\end{equation}
If $\mathcal{F}_\bbR=\LpR$, $1\le p<\infty$, or
$\mathcal{F}_\bbR=C_0(\bbR,X)$, the
space of continuous functions vanishing at infinities (or another
Banach function space as in~\cite{RS1}) this defines a strongly
continuous  semigroup of operators $\{E_\bbR^t\}_{t\ge0}$ whose
generator will be  denoted by $\Gamma_\bbR$.  In the case
$\mathcal{F}_\bbR=\LpR$, $\Gamma_\bbR$ is the closure (in
$\LpR$) of the operator $-d/dt +\calA$ where $(\calA f)(t)= Af(t)$ and
 \begin{eqnarray*}
 \Dom(-d/dt&+&\calA)=\Dom(-d/dt)\cap\Dom (\calA)\\&=&
\{v\in \LpR : v\in AC(\bbR,X), \, v'\in \LpR,\\
&&\quad \text{$v(s)\in \Dom(\calA)$
for almost every~$s$, and $-v' +Av \in \LpR$}  \}.
\end{eqnarray*}
The important properties of this ``evolution semigroup"
 are summarized in the following remarks;
 see \cite{LMS2} and also further developments in
\cite{vanNbook,RS1,RolandDis}.
The unit circle in $\bbC$ is denoted here by $\bbT=\{z\in\bbC:|z|=1\}$.
\begin{rmks}\lb{AUTsmt}
The spectrum $\sigma(E_\bbR^t)$, for $t>0$, is invariant with respect to
rotations centered at the origin, and $\sigma(\Gamma_\bbR)$ is invariant
with respect to translations along $i\bbR$.  Moreover, the following are
equivalent.
\begin{enumerate}
\item $\sigma(e^{tA})\cap\bbT=\emptyset$ on $X$;
\item $\sigma(E_\bbR^t)\cap\bbT=\emptyset$ on $\mathcal{F}_\bbR$;
\item $0\in\rho(\Gamma_\bbR)$ on $\mathcal{F}_\bbR$.
\end{enumerate}
As a consequence,
\begin{equation}\label{ATsmt}
\sigma(E_\bbR^t)\setminus\{0\}=e^{t\sigma(\Gamma_\bbR)},\quad t>0.
\end{equation}
\end{rmks}
\noindent Note that $\{E_\bbR^t\}_{t\ge0}$ has the spectral mapping
property~\eqref{ATsmt} on
$\mathcal{F}_\bbR$ even if the underlying semigroup
$\{e^{tA}\}_{t\ge0}$ does not have the spectral mapping property on
$X$.  In the latter case, it may be that the  exponential stability of
the solutions to $\dot x=Ax$ on $\bbR$ are not determined by the
spectrum of $A$.  However,  such stability {\em is determined} by the
spectrum of $\Gamma_\bbR$.  This is made explicit by the following
corollary of Remarks~\ref{AUTsmt}: The spectral bound
$s(\Gamma_\bbR)$ and the growth bound
$\omega_0(E_\bbR^t)$
for the evolution semigroup coincide and are equal to the growth bound
of   $\{e^{tA}\}_{t\ge0}$\,:
\begin{equation*}\label{GSbounds}
s(\Gamma_\bbR)=\omega_0(E_\bbR^t)=\omega_0(e^{tA}).
\end{equation*}

One of the difficulties related to nonautonomous problems is that their
associated evolution families are two-parameter families of
operators. From this point of view, it would be of interest to define a
{\em
one-parameter} semigroup that is associated to the solutions of the
nonautonomous Cauchy problem \eqref{nacp}.  For such a semigroup to be
useful, its properties should be closely connected to the asymptotic
behavior of the original nonautonomous problem.  Ideally, this semigroup
would have a generator that plays the same significant role in
determining the stability of the solutions as the operator $A$ played in
Lyapunov's classical stability theorem for finite-dimensional autonomous
systems, $\dot x=Ax$.  This can, in fact, be done and the operator of
interest is the generator of the following {\em evolution semigroup}
that is induced by the two-parameter evolution
family:  if $\U=\{U(t,\tau)\}_{t\ge\tau}$ is an
evolution family, define operators $E_\bbR^t$, $t\ge0$, on
 $\mathcal{F}_\bbR=\LpR$ or
$\mathcal{F}_\bbR=C_0(\bbR,X)$ by
\begin{equation}\label{TVevolsgR}
(E^t_\bbR f)(\tau) = U (\tau, \tau -t)f(\tau -t), \quad \tau\in \bbR,
\quad t\ge 0.
\end{equation}
When $\mathcal{U}$  exponentially bounded,
this defines a strongly continuous evolution semigroup on
$\mathcal{F}_\bbR$ whose generator will be denoted by $\Gamma_\bbR$.
As shown in \cite{LMS2} and \cite{RS1} the spectral mapping theorem
\eqref{ATsmt} holds for this semigroup.
Moreover, the existence of an exponential dichotomy for solutions to
$\dot{x}(t) = A(t)x(t),\quad t\in \bbR$,
is characterized by the condition that $\Gamma_\bbR$ is
invertible on $\mathcal{F}_\bbR$.   Note that in the autonomous case
where $U(t,\tau)=e^{(t-\tau)A}$, this is the evolution semigroup defined
in
\eqref{AUTevolsgR}.   In the nonautonomous case, the construction of
an evolution semigroup is a way to
 ``autonomize'' a time-varying Cauchy problem by replacing the
time-dependent differential equation $\dot{x}=A(t)x$ on $X$ by an
autonomous differential equation $\dot{f}=\Gamma f$ on a
superspace of $X$-valued functions.

%---------------------- section 2 -----------------------
\section{Evolution Semigroups and Cauchy Problems}

In order to tackle the problem of characterizing the exponential
stability of solutions to the nonautonomous Cauchy problem
\eqref{nacp} on the half-line, $\bbR_+$,
the following variant of the above evolution semigroup is needed.
As before, let
$\{U(t,\tau)\}_{t\ge\tau}$ be an exponentially bounded
evolution family, and define operators $E^t$, $t\ge0$, on
functions $f:\bbR_+\to X$  by
\begin{equation}\lb{TVevolsgHL}
(E^tf)(\tau)= \begin{cases}
  U(\tau,\tau-t)f(\tau-t),\quad &0 \le t \le \tau \\
0,\qquad &0 \le \tau < t.
  \end{cases}
\end{equation}
This defines a strongly continuous semigroup of operators on
the space of functions $\mathcal{F}=\LpRp$, and the generator
of this {\em evolution semigroup} will be  denoted by
$\Gamma$.   This also defines a strongly continuous semigroup on
$C_{00}(\bbR_+,X)=\{f\in C_0(\bbR_+,X):f(0)=0\}$.
For more information on evolution semigroups on the half-line
 see also \cite{MiRaSc,vanN1,vanNbook,RolandDis}.

\subsection{Stability}
The primary goal of this subsection is
to identify the useful properties of the semigroup of operators
defined in \eqref{TVevolsgHL}
which will be used in the subsequent sections.  In particular,
the following spectral mapping theorem will allow this semigroup
to be used in characterizing the exponential stability of
solutions to~\eqref{nacp} on $\bbR_+$.
See also \cite{RS1,RolandDis} for different proofs.
The spectral symmetry portion of this theorem is due to
R.~Rau \cite{Rau}.

\begin{thm}\lb{SMT+}
Let $\mathcal{F}$ denote $C_{00}(\bbR_+,X)$ or $\LpRp$.
The spectrum $\sigma(\Gamma)$ is a half plane, the spectrum
$\sigma(E^t)$ is a disk centered at the origin, and
\begin{equation}\lb{SMT+eq}
e^{t\sigma(\Gamma)}=\sigma(E^t)\setminus\{0\},\quad t>0.
\end{equation}
\end{thm}

\begin{proof}
The arguments for the two cases $\mathcal{F}=\Coo$ and
$\mathcal{F}=\LpRp$ are similar, so only the first one
is considered here.

We first note that $\sigma(\Gamma)$ is invariant under translations
along $i\bbR$ and $\sigma(E^t)$ is invariant under rotations about
zero.  This spectral symmetry is a consequence of
the fact that  for $\xi\in\bbR$,
\begin{equation}\lb{rescale}
E^t e^{i\xi\cdot}f=e^{i\xi\cdot} e^{-i\xi t}E^t f,
\quad\mbox{and}\quad
\Gamma e^{i\xi\cdot}=e^{i\xi\cdot}(\Gamma-i\xi).
\end{equation}

The inclusion $e^{t\sigma(\Gamma)}\subseteq
\sigma(E^t)\setminus\{0\}$ follows from the standard spectral
inclusion for strongly continuous semigroups \cite{Nagel}.
In view of the spectral symmetry, it suffices to
show that $\sigma(E^t)\cap\bbT=\emptyset$ whenever $0\in\rho(\Gamma)$.
To this end, we replace the Banach space $X$ in
Remarks~\ref{AUTsmt}  by $\Coo$ and consider two semigroups
$\{\tilde{E}^t\}_{t\geq 0}$ and
$\{\mathcal{E}^t\}_{t\geq 0}$ with
generators $\tilde{\Gamma}$ and $\mathcal{G}$,
respectively, acting on the space
$C_0(\bbR,C_{00}(\bbR_+,X))$.  These semigroups are defined by
\begin{eqnarray*}
(\tilde{E}^th)(\tau,\theta)&=& \begin{cases} U(\theta,
\theta-t)h(\tau -t,\theta-t)&\text{for}\quad
\theta \geq t,\\
0 &\text{for}\quad 0\le \theta <t,
\end{cases}\\
(\mathcal{E}^th)(\tau ,\theta)&=&
\begin{cases}
U(\theta ,\theta -t)h(\tau ,\theta -t)
&\text{for}\quad \theta \geq t,
\\ 0 , & \text{for}\quad 0\leq
\theta\le t,
\end{cases}
\end{eqnarray*}
where $\tau \in \bbR$ and $h(\tau ,\cdot )\in C_{00}(\bbR_+, X)$.
Note that if $H\in C_0(\bbR,C_{00}(\bbR_+,X))$, then
$h(\tau,\cdot):=H(\tau)\in C_{00}(\bbR_+,X)$ and we recognize
$\{\tilde{E}^t\}_{t\geq 0}$ as the evolution semigroup induced by
$\{E^t\}_{t\geq 0}$, as in \eqref{AUTevolsgR}:
$$
(\tilde{E}^tH)(\tau )=E^t H(\tau -t).
$$
Also, the semigroup $\{\mathcal{E}^t\}_{t\ge0}$ is the family of
multiplication operators given by
$$
(\mathcal{E} ^tH)(\tau)=E^t H(\tau).
$$
The generator $\mathcal{G}$ of this semigroup is the operator of
multiplication
by $\Gamma$: \
$(\mathcal{G}H)(\tau )=\Gamma(H(\tau ))$,
where $H(\tau )\in \Dom(\Gamma)$ for $\tau \in \bbR$.
In particular, if $0\in \rho(\Gamma)$ on $\mathcal{F}$, then
$(\mathcal{G}^{-1}H)(\tau )=\Gamma^{-1}(H(\tau ))$, and so
$0\in \rho(\mathcal{G})$.

Let $J$ denote the isometry on $C_0(\bbR,C_{00}(\bbR_+,X))$
given by $(J h)(\tau ,\theta )=h(\tau+\theta ,\theta )$ for $\tau \in
\bbR$, $\theta\in\bbR_+$.  Then $J$ satisfies the identity:
$$
(\mathcal{E}^tJh)(\tau,\theta)=
(J\tilde{E}^th)(\tau,\theta),\quad \tau \in \bbR,\quad
\theta \in \bbR_+.
$$
It follows that
$\mathcal{G}JH=J\tilde{\Gamma }H$ for
$H\in \Dom(\tilde{\Gamma })$, and
$J^{-1}\mathcal{G}H=\tilde{\Gamma }J^{-1}H$ for
$H\in \Dom(\mathcal{G})$.  Consequently
$\sigma (\mathcal{G})=\sigma (\tilde{\Gamma })$
on $C_0(\bbR,C_{00}(\bbR_+,X))$. In particular, $0\in \rho
(\tilde{\Gamma })$.  Therefore,
$\sigma (E^t)\cap \bbT=\emptyset$ follows from Remarks~\ref{AUTsmt}
applied to the semigroup $\{E^t\}_{t\ge0}$ on $\mathcal{F}$ in place of
$\{e^{tA}\}_{t\ge0}$ on $X$.

The facts that  $\sigma(\Gamma)$ is a half plane and
$\sigma(E^t)$ is a disk follow from the spectral mapping property
\eqref{SMT+eq} and \cite[Proposition 2]{Rau}.
\end{proof}

An important consequence of this theorem is the property that the growth
bound
$\omega_0(E^t)$ equals the spectral bound $s(\Gamma)$.  This leads to
the following simple result on stability.

\begin{thm}\lb{ExpStabThm+} Let $\mathcal{F}$ denote $C_{00}(\bbR_+,X)$
or $\LpRp$.  An exponentially bounded evolution family
$\{U(t,\tau)\}_{t\ge\tau}$ is exponentially stable if and only if the
growth bound
$\omega_0(E^t)$ of the induced evolution semigroup on $\mathcal{F}$
is negative.
\end{thm}

\begin{proof}  Let $\mathcal{F}=\Coo$.  If $\{U(t,\tau)\}_{t\ge\tau}$ is
exponentially stable, then there exist $M>1$,
$\beta>0$ such that $\|U(t,\tau)\|_{\mathcal{L}(X)}
\le Me^{-\beta(t-\tau)}$, $t\ge
\tau$. For $\t\ge0$ and $f\in\Coo$,
$$
\begin{aligned} \|E^\tau f\|_{\Coo} = \sup_{t>0}\|E^\tau f(t)\|_X
&=\sup_{t>\tau}\|U(t,t-\tau)f(t-\t)\|_X \\
&\le\sup_{t>\tau}\|U(t,t-\tau)\|_{\mathcal{L}(X)} \|f(t-\tau)\|_X \\
&\le Me^{-\beta \tau} \|f\|_{\Coo}.
\end{aligned}
$$

Conversely, assume there exist $M>1$, $\alpha>0$ such that
$\|E^t\|\le Me^{-\alpha t}$, $t\ge 0$.  Let $x\in X$,
$\|x\|=1$.  For fixed $t> \tau > 0$, choose $f\in\Coo$ such that
$\|f\|_{\Coo}=1$ and $f(\tau)=x$.  Then,
$$
\begin{aligned}
\|U(t,\tau )x\|_X=\|U(t,\tau )f(\tau )\|_X &=\|E^{(t-\tau )}f(t)\|_X \\
&\le\sup_{\theta>0}\|E^{(t-\tau)}f(\theta)\|_X\\
&=\|E^{(t-\tau)}f\|_{\Coo}\\
&\le Me^{-\alpha(t-\tau)}.
\end{aligned}
$$

A similar argument works for $\mathcal{F}=\LpRp$.
\end{proof}

The remainder of this subsection focuses on the operator
used for determining exponential stability.  In fact, stability
is characterized by the boundedness of this operator which,
as seen below, is equivalent to the invertibility of
$\Gamma$, the generator of the evolution semigroup.
We begin with the autonomous case.

R. Datko and J.~van Neerven have characterized the exponential
stability of solutions for autonomous equations $\dot x=Ax$,
$t\ge0$, in terms of a convolution operator,
${\mathbb G}$, induced by $\{e^{tA}\}_{t\ge0}$.
In this autonomous setting,
\begin{equation}\lb{AUTevolsgHL}
(E^t f)(\tau)= \begin{cases}
  e^{t A} f(\tau-t),\quad &0 \le t \le \tau \\
0,\qquad &0 \le \tau < t,
  \end{cases}
\end{equation}
and the convolution operator takes the following form: for
$f\in L^1_{loc}(\bbR_+,X)$,
\begin{equation}\lb{AUTbbG}
(\Bbb{G}f)(t):=\int_0^t e^{\tau A}f(t-\tau)\,d\tau
=\int_0^\infty (E^\tau f)(t)\,d\tau ,\quad t\ge 0.
\end{equation}
For reader's convenience we cite Theorem~1.3 of \cite{vanN1}
(see also \cite{Datko}) in the following remarks.

\begin{rmks}\label{AUTvanN} If $\{e^{tA}\}_{t\ge0}$ is a strongly
continuous semigroup on $X$, and $1\le p<\infty$,  then the
following are equivalent:
\begin{enumerate}
\item $\omega_0(e^{tA})<0$;
\item $\Bbb{G}f\in\LpRp$ for all $f\in\LpRp$;
\item $\Bbb{G}f\in\CoRp$ for all $f\in\CoRp$.
%\item $\displaystyle{\sup_{s>0}
%      \left\|\int_0^s e^{tA}f(t)\,dt\right\|_X <\infty}$
%       for all $f\in\CoRp$.
\end{enumerate}
\end{rmks}

\begin{rmks}\label{vNTR}
\begin{enumerate}
\item[(a)]  Note  that condition (ii)
is equivalent to the boundedness
of $\Bbb{G}$ on $L^p(\bbR_+,X)$. To see this, it
suffices to show that the map $f\mapsto \Bbb{G}f$ is a
closed operator on $L^p(\bbR_+,X)$, and then apply the
closed graph theorem.  For this, let $f_n\to f$ and
$\Bbb{G}f_n\to g$ in $L^p(\bbR_+,X)$.  Then
$(\Bbb{G}f_n)(t)\to (\Bbb{G}f)(t)$ for each $t\in \bbR$.
Also,  every norm-convergent sequence in $L^p(\bbR_+,X)$
contains a subsequence that converges pointwise almost
everywhere. Thus, $(\Bbb{G}f_{n_k})(t)\to g(t)$ for almost
all $t$.  This implies that $\Bbb{G}f=g$, as claimed.

\item[(b)]  Also, condition~(iii) is
equivalent to the boundedness of $\Bbb{G}$
on $C_0(\bbR_+,X)$. This follows from
the uniform boundedness principle applied
to the  operators $\bbG_t :f\mapsto
\int_0^te^{\tau A}f(t -\tau)\,d\tau$.
\end{enumerate}
\end{rmks}

We now extend this result so that it may be used to
describe exponential stability for a nonautonomous equation.
For this define
an operator $\bbG$ in an analogous way: let
$\{U(t,\tau)\}_{t\ge\tau}$ be an evolution family  and
$\{E^t\}_{t\ge0}$ the evolution semigroup in \eqref{TVevolsgHL}.
Then define $\bbG$ for $f\in L^1_{loc}(\bbR_+,X)$ as
\begin{equation}\label{TVbbG}
\begin{aligned}
(\bbG f)(t):=&\int_0^\infty (E^\tau f)(t)\,d\tau
=\int_0^t U(t,t-\tau)f(t-\tau)\,d\tau  \\
=&\int_0^t U(t,\tau)f(\t)\,d\tau,\qquad t\ge0.
\end{aligned}
\end{equation}
For $\bbG$ acting on $\mathcal{F}=C_{00}(\bbR_+,X)$ or $\LpRp$,
standard semigroup properties show that $\bbG$ equals
$\,-\Gamma^{-1}$ provided the semigroup $\{E^t\}_{t\ge 0}$ or
the evolution family is uniformly exponentially stable.
Parts (i) $\Leftrightarrow$ (ii) of Remarks~\ref{AUTvanN}
 and the nonautonomous version below are the classical results
by R.~Datko \cite{Datko}.  Our proof uses the evolution semigroup and
creates
a formally autonomous problem so that Remarks~\ref{AUTvanN} can be
applied.

\begin{thm}\label{TVvanN} The following are equivalent for the
evolution family of operators $\{U(t,\tau)\}_{t\ge\tau}$ on $X$.
\begin{enumerate}
\item $\{U(t,\tau)\}_{t\ge\tau}$ is exponentially stable;
\item $\Bbb{G}$ is a bounded operator on $\LpRp$;
\item $\Bbb{G}$ is a bounded operator on $\CoRp$.
\end{enumerate}
\end{thm}

\noindent Before proceeding with the proof,
note that statement~(ii) is equivalent
to the statement: $\Bbb{G}f\in
L^p(\bbR_+,X)$ for each $f\in
L^p(\bbR_+,X)$. This is seen as
in Remark~\ref{vNTR}, above. See also \cite{Buse}
for similar facts.

\begin{proof}
By Theorem~\ref{ExpStabThm+}, (i) implies
that $\{E^t\}_{t\ge0}$ is exponentially stable, and
formula \eqref{TVbbG} implies (ii)
and (iii).  The implication (ii)$\Rightarrow$(i)
will be proved here; the argument for (iii)$\Rightarrow$(i) is
similar.  The main idea is again to use the ``change-of-variables''
technique, as in the proof of
Theorem~\ref{SMT+}.

Consider the operator $\tilde{\Bbb{G}}$
on $L^p(\bbR,L^p(\bbR_+,X)) =L^p(\bbR\times \bbR_+,X)$
defined as multiplication by $\Bbb{G}$.  More precisely,
for $h\in L^p(\bbR\times \bbR_+,X)$ with
${\bf h}(\theta):=h(\theta,\cdot)\in L^p(\bbR_+,X)$, define
$$
(\tilde{\Bbb{G}}h)(\theta ,t )=\Bbb{G}({\bf h}(\theta))(t )=
\int^t _0U(t  ,t -\tau)h(\theta ,t  -\tau)\,d\tau,\quad
t  \in \bbR_+,\quad \theta \in\bbR.
$$

In view of statement (ii), this operator is bounded.
For the isometry $J$ defined on the space
$L^p(\bbR,L^p(\bbR_+,X))$ by
$(Jh)(\theta ,t  )=h(\theta +t ,t  )$, we have
\begin{equation}\label{JinvGJ}
(J^{-1} \tilde{\Bbb{G}}Jh)(\theta,t )
=\int^t _0 U(t  ,t  -\tau)
h(\theta -\tau,t  -\tau)\,d\tau.
\end{equation}

Next, let  $\{E^t\}_{t\ge0}$ be the evolution
semigroup~\eqref{TVevolsgHL} induced by
$\{U(t,\tau)\}_{t\ge\tau}$, and define $\bbG_*$ to be the operator
of convolution with this semigroup as in \eqref{AUTbbG}; that is,
\begin{equation}\label{bbG*}
(\bbG_*{\bf h})(\theta)=
\int^\infty_0 E^\tau {\bf h}(\theta -\tau)\,d\tau,\quad
{\bf h}\in L^p(\bbR,L^p(\bbR_+,X)).
\end{equation}
If $h(\theta,\cdot)={\bf h}(\theta)\in L^p(\bbR_+,X)$, then by
definition~\eqref{TVevolsgHL}, evaluating \eqref{bbG*} at $t$ gives
\begin{equation}\label{bbG*2}
\left[(\bbG_*{\bf h})(\theta)\right](t)=(\bbG_*h)(\theta ,t  )=
\int^t _0U(t ,t -\tau)h(\theta -\tau,t  -\tau)\,d\tau,
\quad t  \in \bbR_+,\ \theta \in \bbR.
\end{equation} From \eqref{JinvGJ} it follows that
$\bbG_*=J^{-1}\tilde{\Bbb{G}}J$
is a bounded operator on $L^p(\bbR,L^p(\bbR_+,X))$.

Now, each function ${\bf h}_+\in L^p(\bbR_+,L^p(\bbR_+,X))$ is
an $L^p(\bbR_+,X)$-valued function
on the half line $\bbR_+$. We extend each such ${\bf h}_+$ to a
function ${\bf h}\in L^p(\bbR,L^p(\bbR_+,X))$ by setting
${\bf h}(\theta)={\bf h}_+(\theta )$ for $\theta \geq 0$ and
${\bf h}(\theta)=0$ for $\theta <0$.
Note that $\bbG_*{\bf h}\in L^p(\bbR,L^p(\bbR_+,X))$ because
$\bbG_*$ is bounded on
$L^p(\bbR,L^p(\bbR_+,X))$.
Consider the function
${\bf f}_+:\bbR_+\to L^p(\bbR_+,X)$ defined by
$$
{\bf f}_+(t)=\int_0^t E^\tau{\bf h}_+(t
-\tau)\,d\tau=\int_0^\infty E^\tau{\bf h}(t-\tau)\,d\tau,
\quad t \in \bbR_+.
$$

\noindent To complete the proof of the theorem, it suffices to prove
the following claim:  $${\bf f}_+\in L^p(\bbR_+,L^p(\bbR_+,X)).$$
Indeed, the operator ${\bf h}_+\mapsto {\bf f}_+$
is the convolution operator as in
\eqref{AUTbbG} defined by the semigroup operators
$E^t$ instead of $e^{tA}$.
An application of Remarks~\ref{AUTvanN} to
$E^t$ on $L^p(\bbR_+,X)$ (in place of $e^{tA}$ on $X$)
shows that
the semigroup $\{E^t\}_{t\ge0}$ is exponentially stable
on $L^p(\bbR_+,X)$ provided:
\[
{\bf f}_+\in L^p(\bbR_+,L^p(\bbR_+,X)) \quad \text{for
each}\ {\bf h}_+\in L^p(\bbR_+,L^p(\bbR_+,X)).
\]
\noindent But if $\{E^t\}_{t\ge0}$ is exponentially stable,
the evolution family $\{U(t,\tau)\}_{t\ge\tau}$ is
exponentially stable by Theorem~\ref{ExpStabThm+}.

To prove the claim, apply formula \eqref{bbG*2} for
$h(\theta ,t )=h_+(\theta ,t  ),\theta \geq 0$ and
$h(\theta ,t  )=0$, $\theta <0$,
$t  \in \bbR_+$, where
${\bf h}_+(\theta )=h_+(\theta ,\cdot )$. This gives
$$
(\bbG_*h)(\theta ,t )=
\begin{cases}
\int^{\min \{\theta ,t \}}_0
U(t ,t -\tau)h_+(\theta-\tau,t  -\tau)\,d\tau \quad
&\text{for}\quad \theta\ge 0,\  t  \in \bbR_+\\
(\bbG_*h)(\theta ,t )=0,\quad
&\text{for}\quad \theta <0, \  t  \in \bbR_+
\end{cases}
$$
Thus, the function
$$
\theta \mapsto (\bbG_*h)(\theta ,\cdot)
   =(\bbG_*{\bf h})(\theta )\in L^p(\bbR_+,X)
$$
is in the space
$L^p(\bbR_+,L^p(\bbR_+,X))$.
On the other hand, denoting $f_+(\theta ,\cdot):={\bf f}_+(\theta)
\in L^p(\bbR_+,X)$, we have that
$$
f_+(\theta ,t  )=\int^{\min \{\theta ,t  \}}_0
U(t,t-\tau)h_+(\theta -\tau,t-\tau)\,d\tau,\quad
\theta ,t  \in \bbR_+.
$$
Thus, $\theta \mapsto f_+(\theta ,\cdot)
=(\bbG_*h)(\theta ,\cdot )$
is a function in $L^p(\bbR_+,L^p(\bbR_+,X))$,
and the claim is proved.
\end{proof}

This theorem makes explicit, in the case of
the half line $\bbR_+$, the relationship between the stability of an
evolution family $\{U(t,\tau)\}_{t\ge\tau}$ and the generator, $\Gamma$,
of the corresponding evolution semigroup \eqref{TVevolsgHL}.  Indeed, as
shown above, stability is equivalent to the boundedness of $\Bbb{G}$, in
which case $\Bbb{G}=-\Gamma^{-1}$.   Combining Theorems \ref{SMT+},
\ref{ExpStabThm+} and \ref{TVvanN}
yields the following corollary.

\begin{cor}\lb{vanNCor}
Let $\{U(t,\tau)\}_{t\ge\tau}$ be an exponentially
bounded evolution family and let $\Gamma$ denote the
generator of the induced evolution semigroup on
$L^p(\bbR_+,X)$, $1\le p<\infty$, or $C_{00}(\bbR_+,X)$.
The following are equivalent:
\begin{enumerate}
\item $\{U(t,\tau)\}_{t\ge\tau}$ is exponentially stable;
\item $\Gamma$ is invertible with $\Gamma^{-1}=-\Bbb{G}$;
\item $s(\Gamma)<0$.
\end{enumerate}
\end{cor}
\noindent For more information on stability and dichotomy of
evolution families on the semiaxis see \cite{MiRaSc}.

\subsection{Perturbations and robust stability}%-----------------

This subsection briefly considers
perturbations of \eqref{nacp} of the form
\begin{equation}\lb{nacp+D}
\dot x(t)=(A(t)+D(t))x(t), \quad t\ge0.
\end{equation}
It will not, however, be assumed that \eqref{nacp+D} has a
differentiable solution.  For example,
 let $\{e^{tA_0}\}_{t\ge0}$ be a strongly continuous semigroup
generated by $A_0$, let $A_1(t)\in \mathcal{L}(X)$ for $t\ge0$,
and define $A(t)=A_0+A_1(t)$.  Then even if $t\mapsto A_1(t)$
is continuous, the Cauchy problem \eqref{nacp} may   not have
a differentiable solution for all initial conditions $x(0)=x\in
\Dom(A)=\Dom(A_0)$ (see, e.g.,\cite{Phillips}).  Therefore we will
want our development to allow for equations with solutions that
exist only in the following mild sense.

Let $\{U(t,\tau)\}_{t\ge\tau}$ be an
evolution family of operators corresponding to a solution  of
\eqref{nacp}, and consider the nonautonomous inhomogeneous equation
\begin{equation}\lb{inhomog}
\dot x(t)=A(t)x(t)+f(t),\qquad t\ge0,
\end{equation}
where $f$ is a locally integrable $X$-valued function on $\bbR_+$.
A function $x(\cdot)$ is a mild solution of \eqref{inhomog}
with initial value $x(\theta)=x_\theta \in \Dom(A(\theta))$ if
\[
x(t)=U(t,\theta)x_\theta+\int_\theta^tU(t,\tau)f(\tau)\,d\tau,
\quad t\ge \theta.
\]
Given operators $D(t)$, the existence of mild
solutions to an additively perturbed equation \eqref{nacp+D}
corresponds to the existence of an evolution
family  $\{U_1(t,\tau)\}_{t\ge \tau}$ satisfying
\begin{equation}\lb{varpar}
U_1(t,\theta)x=U(t,\theta)x+\int_\theta^tU(t,\t)D(\t)U_1(\t,\theta)x\,d\t.
\end{equation}
for all $x\in X$.
It will be assumed that the perturbation operators, $D(t)$,  are
strongly measurable and essentially bounded functions of
$t$.  In view of this, we use the notation  $\mathcal{L}_s(X)$
to denote the set $\mathcal{L}(X)$  endowed with the
strong operator topology and  use
$L^\infty(\bbR_+,\mathcal{L}_s(X))$ to denote the set of
 bounded, strongly measurable
$\mathcal{L}(X)$-valued functions on $\bbR_+$.
A function $D(\cdot)\in L^\infty(\bbR_+,\mathcal{L}_s(X))$
induces a  multiplication operator $\mathcal{D}$ defined
by $\mathcal{D} x(t)=D(t)x(t)$, for $x(\cdot)\in L^p(\bbR_+,X)$.
In fact,  $\mathcal{D}$ is a bounded operator
on $\LpRX$  with $\|\mathcal{D}\| \le
\|D(\cdot)\|_\infty:=\esssup_{t\in\Bbb{R}_+}\|D(t)\|$.
%%\cite{RolandDis}.  %%[p.~36]

Evolution semigroups induced by an evolution family as in
\eqref{TVevolsgHL} have been studied by several authors
who have characterized such semigroups in
terms of their generators on general Banach function spaces  of
$X$-valued functions (see \cite{RRSV,RolandDis}
and the bibliography therein).  The sets  $\mathcal{F}=\LpRp$ or
$\mathcal{F}=\Coo$  considered here are examples of more general
``Banach
function spaces."   In the development that follows we use
a theorem of  R. Schnaubelt~\cite{RolandDis} (see also R\"abiger et
al.~\cite{RRS, RRSV}) which shows exactly when a strongly continuous
semigroup on $\mathcal{F}$ arises  from a strongly continuous evolution
family on
$X$.  We state a version of this result which will be used below; a
more general version is proven in
\cite{RRSV}.  The set $C_c^1(\bbR_+)$ consists of differentiable
functions on $\bbR_+$ that have compact support.

\begin{thm}\lb{RRS-evsgThm}
Let $\{T^t\}_{t\ge0}$ be a strongly continuous semigroup generated by
$\Gamma$ on $\mathcal{F}$.  The
following are equivalent:
\begin{enumerate}
\item $\{T^t\}_{t\ge0}$ is an evolution semigroup;  i.e.,
there exists an exponentially bounded evolution family so that $T^t$ is
defined as in \eqref{TVevolsgHL};
\item there exists a core, $\scrC$, of $\Gamma$ such that for all
$\vphi\in C_c^1(\bbR_+)$, and $f\in\scrC$, it follows that $\vphi f\in
\Dom(\Gamma)$ and $\Gamma(\vphi f)=-\vphi'f+\vphi\Gamma f$.  Moreover,
there exists $\lambda\in\rho(\Gamma)$ such that
$R(\lambda,\Gamma): \mathcal{F}\to \Coo$ is continuous with dense range.
\end{enumerate}
\end{thm}

Now let $\{U(t,\tau)\}_{t\ge \tau}$ be an evolution family on $X$ and
let
$\Gamma$ be the generator of the corresponding evolution semigroup,
$\{E^t\}_{t\ge0}$, as in \eqref{TVevolsgHL}.
If $D(\cdot)\in L^\infty(\bbR_+,\mathcal{L}_s(X))$, then the
multiplication operator $\mathcal{D}$ is a bounded operator
 on $\mathcal{F}=L^p(\bbR_+,X)$. Since a bounded perturbation of a
generator of a strongly continuous semigroup is itself such a generator,
the operator $\Gamma_1=\Gamma+\mathcal{D}$ generates a strongly
continuous semigroup, $\{E_1^t\}_{t\ge0}$ on
$\mathcal{F}$ (see, e.g., \cite{Pazy}). In fact, $\Gamma_1$ generates
an {\em evolution} semigroup, see \cite{RRS,RolandDis}:

\begin{prop}\lb{GammaRobustStab} Let $D(\cdot)\in
L^\infty(\bbR_+,\mathcal{L}_s(X)).$, and let $\{U(t,\tau)\}_{t\ge \tau}$
be an
exponentially bounded evolution family.   Then there exists a unique
evolution family $\Uone=\{U_1(t,\tau)\}_{t\ge
\tau}$ which solves the integral equation \eqref{varpar}.  Moreover,
$\Uone$ is exponentially stable if and only if $\Gamma+\mathcal{D}$ is
invertible.
\end{prop}

\begin{proof} As already observed,
$\Gamma_1=\Gamma+\mathcal{D}$ generates a strongly continuous semigroup,
$\{E_1^t\}_{t\ge0}$ on ${\mathcal F}$.  To see that this is, in fact, an
evolution semigroup, note that for
$\lambda\in\rho(\Gamma)\cap\rho(\Gamma_1)$,
\[\Range(R(\lambda,\Gamma))=\Dom(\Gamma)
=\Dom(\Gamma+\mathcal{D})
=\Range(R(\lambda,\Gamma_1))\] is dense in $\Coo$.
Also, if $\scrC$
is a core for $\Gamma$, then it is a core for $\Gamma_1$, and so
for $\vphi\in C_c^1(\bbR)$, $f\in\scrC$,
\[
\Gamma_1(\vphi f)=\Gamma(\vphi f)+\mathcal{D}(\vphi f)=
-\vphi'f+\vphi\Gamma f+\vphi \mathcal{D} f
=-\vphi'f+\vphi(\Gamma+\mathcal{D})f.
\]
Consequently, Corollary~\ref{RRS-evsgThm} shows that $\{E_1^t\}_{t\ge0}$
corresponds to an evolutionary family, $\{U_1(t,\tau)\}_{t\ge \tau}$.
Moreover,  $x(t)=U_1(t,\tau)x(\tau)$ is
seen to define a mild solution to \eqref{nacp+D}.  Indeed,
\begin{equation}\lb{varparF}
E_1^tf=E^tf+\int_0^tE^{(t-\t)}\mathcal{D} E_1^\t f \,d\t,
\end{equation}
holds for all $f\in F$.  In particular, for $x\in X$, and any
$\vphi\in C_c^1(\bbR)$, setting $f=\vphi\otimes x$ in \eqref{varparF},
where $\vphi\otimes x(t)=\vphi(t)x$,
and using a change of variables leads to
\[
\vphi(\theta)U_1(t,\theta)x=\vphi(\theta)U(t,\theta)x+
 \vphi(\theta)\int_\theta^tU(t,\t)D(\t)U_1(\t,\theta)x\,d\t.
\]
Therefore, \eqref{varpar} holds for all $x\in X$.

Finally, Theorem~\ref{vanNCor} shows that $\Uone$ is exponentially
stable if and only if $\Gamma_1$ is invertible.
\end{proof}

The existence of mild solutions under bounded perturbations of
this type is well known (see, e.g., \cite{CP}), but an immediate
consequence of the approach given here is the property of
robustness for the stability of $\{U(t,\tau)\}_{t\ge\tau}$.  Indeed, by
continuity  properties of the spectrum of an operator $\Gamma$, there
exists $\epsilon>0$ such that $\Gamma_1$ is invertible whenever
$\|\Gamma_1-\Gamma\|<\epsilon$;  that is, $\{U_1(t,\tau)\}_{t\ge\tau}$
is
exponentially stable whenever $\|D(\cdot)\|_\infty<\epsilon$.
Also, the type of proof presented here can be extended to address
the case of unbounded perturbations.  For an example of this, we refer
to \cite{RRSV}.
Finally, and most important to the present paper, is the fact that this
approach provides insight into the concept of the stability radius.
This topic is studied next.

%---------------------- section 3 -----------------------
\section{Stability Radius}

The goal of this section is to use the previous development to study
the (complex)  stability radius of an exponentially stable system.
Loosely speaking, this is a measurement on the size of the smallest
operator under which the additively perturbed system looses exponential
stability.  This is an important concept for linear systems theory and
was introduced by D.~Hinrichsen and A.~J.~Pritchard  as the basis for a
state-space approach to studying robustness of linear time-invariant
\cite{HP86} and time-varying systems \cite{HIP89,HP94,PT89}. A
systematic study of various stability radii in the spirit of the
current paper has recently be given by
A.~Fischer and J.~van Neerven \cite{FishvN}.

\subsection{General estimates}

In this subsection we give estimates for the stability radius of general
nonautonomous systems on Banach spaces. The perturbations considered
here are additive ``structured" perturbations of output feedback type.
That is, let $U$ and $Y$ be Banach spaces and let $\Delta(t): Y\to U$
denote an unknown disturbance operator.  The operators $B(t):U\to X$ and
$C(t):X\to Y$  describe the structure of the perturbation in the
following (formal) sense:  if $u(t)=\Delta(t)y(t)$ is viewed as a
feedback for the system
\begin{equation}\label{TVcontrolProb}
\begin{aligned}
\dot x(t)&=A(t)x(t)+B(t)u(t),\qquad x(s)=x_s\in \Dom(A(s)),\\
y(t)&=C(t)x(t),\quad t\ge s\ge0,
\end{aligned}
\end{equation}
then the nominal system $\dot x(t)=A(t)x(t)$ is subject to the
structured perturbation:
\begin{equation}\label{TVpert}
\dot x(t)=(A(t)+B(t)\Delta(t)C(t))x(t),\quad t\ge 0.
\end{equation}
In this section $B$ and $C$ do not represent input and output operators,
rather they describe the structure of the uncertainty of the system.
Also, systems considered throughout this paper are not assumed to have
differentiable solutions and so
\eqref{TVpert} is to be interpreted in the mild sense as described in
\eqref{varpar} where $D(t)=B(t)\Delta(t)C(t)$.  Similarly,
\eqref{TVcontrolProb} is interpreted in the mild sense; that is, there
exists a strongly continuous exponentially bounded evolution family
$\{U(t,\tau)\}_{t\ge\tau}$ on a Banach space $X$ which satisfies
\begin{equation}\lb{mildTVsys}
\begin{aligned}
x(t)&=U(t,s)x(s)+
\int_s^tU(t,\tau)B(\tau)u(\tau)\,d\tau,
\\ y(t)&=C(t)x(t), \qquad t\ge s\ge 0.
\end{aligned}
\end{equation}
In the case of  time-invariant systems, equation \eqref{mildTVsys} takes
the form
\begin{equation}\lb{mildAUTsys}
\begin{aligned}
x(t)&=e^{tA}x_0+\int_0^te^{(t-\t)A}Bu(\t)\,d\t,\\
y(t)&=Cx(t), \qquad  t\ge 0,
\end{aligned}
\end{equation}
where $\{e^{tA}\}_{t\ge0}$ is a strongly
continuous semigroup on $X$ generated by $A$, $x(0)=x_0\in\Dom(A)$.

It should be emphasized that we will not address questions concerning
the existence of solutions for a perturbed system \eqref{TVpert} beyond
the point already discussed in Proposition \ref{GammaRobustStab}.  In
view of that proposition, we make the following assumptions:
$B$, $C$ and $\Delta$  are strongly measurable and essentially
bounded functions of $t$; i.e.,
$B(\cdot)\in L^\infty(\bbR_+,\mathcal{L}_s(U,X))$,
$C(\cdot)\in L^\infty(\bbR_+,\mathcal{L}_s(X,Y))$ and
$\Delta(\cdot)\in L^\infty(\bbR_+,\mathcal{L}_s(Y,U))$.  As such, they
induce bounded multiplication operators,
$\mathcal{B}$, $\mathcal{C}$ and $\tilde{\Delta}$ acting on the spaces
$L^p(\bbR_+,U)$, $L^p(\bbR_+,X)$ and $L^p(\bbR_+,Y)$,
respectively.

Next, for an exponentially bounded evolution family
$\{U(t,\tau)\}_{t\ge\tau}$, define the ``input-output" operator
$\bbL$ on functions $u:\bbR_+\to U$ by the rule
\[
(\bbL u)(t)=C(t)\int_0^t U(t,\t)B(\t)u(\t)\, d\t.
\]
Using the above notation note that $\bbL=\mathcal{C}\bbG\mathcal{B}$.
Much of the stability analysis that follows is based on this observation
in combination with Theorem~\ref{TVvanN} which shows that  the operator
$\bbG$ completely characterizes stability of the corresponding
evolution
family.

We now turn to the definition of the stability radius.
For this let  $\U=\{U(t,\tau)\}_{t\ge\tau}$ be an exponentially stable
evolution family on $X$.
Set $\mathcal{D}=\mathcal{B}\tilde{\Delta}\mathcal{C}$ and let
$\U_\Delta=\{U_\Delta(t,\tau)\}_{t\ge\tau}$ denote  the evolution
family corresponding to solutions of the perturbed
equation~\eqref{nacp+D}.  That is, $\U_\Delta$ satisfies
\begin{equation*}
U_\Delta(t,s)x=U(t,s)x+
\int_s^tU(t,\tau)B(\tau)\Delta(\tau)C(\tau)U_\Delta(\tau,s)x\,d\tau,
\qquad x\in X.
\end{equation*}
Define the (complex) {\em stability radius} for
$\U$ with respect to the perturbation structure
$(B(\cdot),C(\cdot))$ as the quantity
$$
\rstab(\U,B,C)=\sup\{r\ge0:\|\Delta(\cdot)\|_\infty\le r\Rightarrow
\U_\Delta\mbox{ is exponentially stable} \}.
$$
This definition applies to both nonautonomous and autonomous systems,
though in the latter case the notation $\rstab(\{e^{tA}\},B,C))$ will be
used to distinguish the case where all the operators except $\Delta(t)$
are independent of $t$.  We will have occasion to consider the
{\em constant stability radius} which is defined for the case in which
$\Delta(t)\equiv \Delta$ is constant; this will be denoted by
$rc_{stab}(\{e^{tA}\},B,C))$
or $rc_{stab}(\U,B,C))$, depending on the context.
The above remarks concerning $\Gamma+\mathcal{D}$ (see
Proposition~\ref{GammaRobustStab}), when combined with
Theorem \ref{ExpStabThm+}, make it clear  that
\begin{equation}\label{rstabGamma}
\rstab(\U,B,C)=\sup\{r\ge0:\|\Delta(\cdot)\|_\infty\le r\Rightarrow
\Gamma+\mathcal{B}\tilde{\Delta}\mathcal{C}\mbox{ is invertible} \}.
\end{equation}

It is well known that for autonomous systems in which
$U$ and $Y$ are Hilbert spaces and $p=2$, the stability radius may
be expressed in terms of the norm of the input-output operator or the
transfer function:
\begin{equation}\lb{rstabHilbeq}
\frac{1}{\|\bbL\|_{\mathcal{L}(L^2)}}=\rstab(\{e^{tA}\},B,C)=
\frac{1}{\sup_{s\in\bbR}\|C(A-is)^{-1}B\|};
\end{equation}
see, e.g.,~\cite[Theorem~3.5]{HP94}.
For nonautonomous  equations, a scalar example given in Example~4.4 of
\cite{HIP89} shows that, in general, a strict inequality
$1/{\|\bbL\|}<\rstab(\U,B,C)$ may hold.
Moreover, even for autonomous systems, when Banach spaces are allowed or
when $p\neq 2$, Example~\ref{RstabIneqEx} and Example~\ref{StephenEx}
below
will show that neither of the equalities in \eqref{rstabHilbeq}
necessarily
hold. Subsection~\ref{AUTsubsection}
below focuses on  autonomous equations and a primary objective
there is to prove the following result.
\begin{thm}\label{autin}
For the general autonomous systems,
\begin{equation}\lb{rstabBanIneq}
\frac{1}{\|\bbL\|_{\mathcal{L}(L^p)}}\le\rstab(\{e^{tA}\},B,C)\le
\frac{1}{\sup_{s\in\bbR}\|C(A-is)^{-1}B\|},
\quad 1\le p<\infty.
\end{equation}
\end{thm}
As seen next, the lower bound here holds for general nonautonomous
systems and may be proven in a very direct way using the make-up of the
operator
$\bbL=\mathcal{C}\Bbb{G}\mathcal{B}$.   This lower bound is also proven
in  \cite[Theorem~3.2]{HP94} using a completely different approach.

\begin{thm}\lb{L0BoundThm}
Assume $\U$ is an exponentially stable evolution family and let $\Gamma$
denote the generator of the corresponding evolution semigroup.  If
$$B(\cdot)\in L^\infty(\bbR_+,\mathcal{L}_s(U,X))\quad\text{and}\quad
C(\cdot)\in L^\infty(\bbR_+,\mathcal{L}_s(X,Y)),$$
then $\bbL$ is a bounded operator from $L^p({\mathbb R}_+,U)$
to $L^p({\mathbb R}_+,Y)$, $1\le p<\infty$, the formula
$${\mathbb L}=\mathcal{C}\Bbb{G}\mathcal{B}
=-\mathcal{C}\Gamma^{-1}\mathcal{B}$$
holds, and
\begin{equation}\lb{LowerRstabBound}
\frac{1}{\|\bbL\|}\le \rstab(\U,B,C).
\end{equation}
In the ``unstructured" case, where $U=Y=X$ and $B=C=I$, one has
\[\bbL=-\Gamma^{-1},\quad\mbox{ and }\quad
\frac{1}{\|\Gamma^{-1}\|} \le \rstab(\U,I,I) \le
\frac{1}{r(\Gamma^{-1})},
\]
where $r(\cdot)$ denotes the spectral radius.
\end{thm}

\begin{proof} Since $\U$ is exponentially stable,
$\Gamma$ is invertible and $\Gamma^{-1}=-\Bbb{G}$.
The required formula for $\bbL$ follows from \eqref{TVbbG}.

Set $\calH:=\Gamma^{-1}\calB\tilde{\Delta}$.  To prove
\eqref{LowerRstabBound},  let $\Delta(\cdot)\in
L^\infty(\bbR_+,\mathcal{L}_s(Y,U))$ and suppose that
$\|\Delta(\cdot)\|_\infty < 1/\|\bbL\|$.  Then
$\|\bbL\tilde{\Delta}\|<1$, and hence
$I-\bbL\tilde{\Delta}=I+\calC\Gamma^{-1}\calB\tilde{\Delta}$ is
invertible on $\LpRY$.  That is,
$I+\calC\calH$ is invertible on $\LpRY$, and hence
$I+\calH\calC$ is invertible on $\LpRX$ (with inverse
$(I-\calH(I+\calC\calH)^{-1}\calC)$).  Now,
\[
\Gamma+\calB\tilde{\Delta}\calC
=\Gamma(I+\Gamma^{-1}\calB\tilde{\Delta}\calC)
=\Gamma(I+\calH\calC)
\]
and so $\Gamma+\calB\tilde{\Delta}\calC$ is
invertible.  It follows from the expression~\eqref{rstabGamma}
that  $1/\|\bbL\|\le \rstab(\U,B,C)$.

For the last assertion, suppose that
$\rstab(\U,I,I)>1/r(\Gamma^{-1})$. Then there exists $\lambda$ such
that $|\lambda|=r(\Gamma^{-1})$ and
$\lambda+\Gamma^{-1}$ is not invertible.
But then setting $\tilde{\Delta}\equiv\frac{1}{\lambda}$ gives
$\|\tilde{\Delta}\|=\frac{1}{|\lambda|}<\rstab(\U,I,I)$, and so
$\Gamma+\tilde{\Delta}=\tilde{\Delta}(\lambda+\Gamma^{-1})\Gamma$ is
invertible, a contradiction.
\end{proof}

\subsection{The transfer function for nonautonomous systems} %-----

In this subsection  we consider a time-varying
version of equation~\eqref{rstabHilbeq} and then observe that
the concept of a transfer function, or frequency-response function,
arises naturally from these ideas.  For this we assume in this
subsection that $X$, $U$ and $Y$ are Hilbert spaces and $p=2$.

Let $\{U(t,\tau)\}_{t\ge\tau}$ be an uniformly exponentially stable
evolution
family and let $\{E^t\}_{t\ge0}$ be the induced evolution semigroup
with generator $\Gamma$ on $L^2(\bbR_+,X)$.
Recall that $\calB$ and $\calC$ denote multiplication operators, with
respective multipliers $B(\cdot)$ and $C(\cdot)$, that act on the
spaces $L^2(\bbR_+,U)$ and $ L^2(\bbR_+,X)$, respectively.
Let $\tilde\calB$ and $\tilde\calC$ denote operators of multiplication
induced by $\calB$ and $\calC$, respectively; e.g.,
$(\tilde\calB \mathbf{u})(t) =\calB(\mathbf{u}(t))$, for
$\mathbf{u}:\bbR_+\to  L^2(\bbR_+,U)$.  Now consider the operator
$\bbG_*$ as defined in equation~\eqref{bbG*} and note that
operator $\bbL_*:=\tilde\calC\bbG_*\tilde\calB$ may be
viewed (formally) as an input-output operator for the ``autonomized"
system: $\dot f=\Gamma f+\calB u,\ \ g=\calC f$, \  where the state
space is $L^2(\bbR_+,X)$. It follows from the
known Hilbert-space  equalities in \eqref{rstabHilbeq} that
\begin{equation*}
\frac{1}{\|\bbL_*\|}=\rstab(\{E^{t}\},\calB,\calC)=
\frac{1}{\sup_{s\in\bbR}\|\mathcal{C}(\Gamma-is)^{-1}\mathcal{B}\|}.
\end{equation*}
Note, however, that the rescaling identities \eqref{rescale} for
$\Gamma$ imply that
\[\|\bbL_*\|=\|\calC(\Gamma-is)^{-1}\calB\|=
\|\calC\Gamma^{-1}\calB\|=\|\bbL\|,\] and so the stability radius for
the evolution semigroup is also $1/\|\bbL\|$.
In view of the above-mentioned nonautonomous scalar example for which
$1/{\|\bbL\|}<\rstab(\U,B,C)$ we see that even though  the evolution
semigroup (or its generator) completely determines the exponential
stability of a system, it  does not provide a formula for the stability
radius.

However, the operator $\mathcal{C}(\Gamma-is)^{-1}\mathcal{B}$
appearing above suggests that the transfer function for
time-varying systems arises  naturally when viewed in the context of
 evolution semigroups.  Several authors have considered
the concept of a tranasfer function for nonautonomous systems
but the work of J.~Ball, I.~Gohberg, and M.A.~Kaashoek~\cite{BGK}
seems to be the most
comprehensive in providing a system-theoretic input-output
interpretation for the value of such a transfer function at a point.
Their interpretation justifies the term   {\em frequency response}
function for  time-varying finite-dimensional systems with
``time-varying complex exponential inputs."  Our remarks concerning the
frequency response for time-varying infinite-dimensional systems will be
restricted to inputs of the form $u(t)=u_0e^{\lambda t}$.

For motivation, consider the input-output operator $\bbL$ associated
with an autonomous system~\eqref{mildAUTsys} where the nominal system is
exponentially stable.  The {\em transfer function} of $\bbL$ is the
unique bounded analytic
$\mathcal{L}(U,Y)$-valued function, $H$, defined on
$\bbC_+=\{\lambda\in{\mathbb C}: \text{Re }\lambda> 0\}$ such that for
any
$u\in L^2(\bbR_+,U)$,
$$
(\widehat{\bbL u})(\lambda)=H(\lambda)\hat{u}(\lambda),\qquad
\lambda\in\bbC_+,
$$
where \, $\widehat{\ }$ \, denotes the
Laplace transform (see, e.g., \cite{Weiss1}).  In this autonomous
setting,
$A$ generates a uniformly exponentially stable strongly continuous
semigroup, and $\bbL=\calC\bbG\calB$ where
$\bbG$ is the operator of convolution with the semigroup
operators $e^{tA}$ (see~\eqref{AUTbbG}).
Standard arguments show that
$(\widehat{\bbL u})(\lambda)=C(\lambda-A)^{-1}B\hat u(\lambda)$; that
is,
$H(\lambda)=C(\lambda-A)^{-1}B$.

Now let $\bbL$ be the input-output operator for the nonautonomous system
\eqref{mildTVsys}.
We wish to identify the transfer function of $\bbL$ as the Laplace
transform of the appropriate operator.  We are guided by the fact
that, just as $(\lambda-A)^{-1}$ may be expressed as the Laplace
transform of
the semigroup generated by $A$, the operator $(\lambda-\Gamma)^{-1}$
is the Laplace transform of the evolution semigroup.  For nonautonomous
systems, $\bbL$ is again given by $\calC\bbG\calB$, although now $\bbG$
from \eqref{TVbbG}  is not, generally, a convolution operator.  So
instead
recall the operator $\bbG_*$ from \eqref{bbG*} which {\em is} the
operator of convolution with the evolution semigroup $\{E^t\}_{t\ge 0}$.
As noted above, the operator $\bbL_*:=\tilde\calC\bbG_*\tilde\calB$ may
be
viewed  as an input-output operator for an autonomous
system (where the state space is $L^2(\bbR_+,X)$).  Therefore, the
autonomous theory applies directly to show that, for $\mathbf{u}\in
L^2(\bbR_+,L^2(\bbR_+,U))$,
\begin{equation}\label{LapL}
(\widehat{\bbL_*\mathbf{u}})(\lambda)=
\calC(\lambda-\Gamma)^{-1}\calB\hat{\mathbf{u}}(\lambda).
\end{equation}
In other words, the transfer function for $\bbL_*$ is
$\calC(\lambda-\Gamma)^{-1}\calB$, where
$$
\calC(\lambda-\Gamma)^{-1}\calB u=\calC\int_0^\infty
e^{-\lambda\tau}E^\tau\calB u\,d\tau,\qquad u\in L^2(\bbR_+,U).
$$
Evaluating these expressions at $t\in\bbR_+$ gives
\begin{equation}\label{nonauttransfcn}
[\calC(\lambda-\Gamma)^{-1}\calB u](t)=\int_0^t
C(t)U(t,\tau)B(\tau)u(\tau)e^{-\lambda(t-\tau)}\,d\tau.
\end{equation}
It is natural to call $\calC(\lambda-\Gamma)^{-1}\calB$ the transfer
function for the nonautonomous system.  Moreover,
the following remarks show that,
by looking at the right-hand side of \eqref{nonauttransfcn},
this gives a natural ``frequency response'' function for
nonautonomous systems.
To see this, we first consider  autonomous systems and
note that the definition of the transfer function for an
autonomous system can be extended to allow for
a class of ``Laplace transformable" functions that are in
$L^2_{loc}(\bbR_+,U)$ (see, e.g.,~\cite{Weiss1}).  This class
includes constant functions of the form $v_0(t)=u_0$, $t\ge0$, for a
given $u_0\in U$.  If a periodic input signal of the form
$u(t)=u_0e^{i\omega t}$, $t\ge0$, (for some $u_0\in U$ and
$\omega\in{\mathbb R}$) is fed into an autonomous system  with
initial condition $x(0)=x_0$, then, by definition of the input-output
operator, we have
\begin{equation*}\label{freq}
(\bbL u)(t)=C(i\omega-A)^{-1}B u_0\cdot e^{i\omega t}-Ce^{tA}x_0,
\quad (\bbL u)(t)=C\int_0^te^{(t-s)A}Bu(s)\,ds.
\end{equation*}
Thus, the output
$$
y(t;u(\cdot),x_0)=(\bbL u)(t)+Ce^{tA}x_0=
C(i\omega-A)^{-1}B u_0\cdot e^{i\omega t}
$$
has the same frequency  as the input.  In view of this, the function
$C(i\omega-A)^{-1}B$ is sometimes called the frequency response
function.  Now recall that the semigroup $\{e^{tA}\}_{t\ge0}$ is
stable and so $\lim_{t\to\infty}\|Ce^{tA}x_0\|=0$.  On the other hand,
consider $v(t)=u_0$ and (formally) apply
$\calC(i\omega-\Gamma)^{-1}\calB$ to this $v$.
For $x_0=(i\omega-A)^{-1}Bu_0$, a calculation based on the Laplace
transform formula for the resolvent of the generator (applied to the
evolution
semigroup $\{E^t\}_{t\ge0}$) yields the identity
$$\left[\calC(i\omega-\Gamma)^{-1}\calB u_0\right](t)=
C(i\omega-A)^{-1}B u_0-Ce^{tA}x_0\cdot e^{-i\omega t}.
$$

Let us consider this expression
$[\calC(i\omega-\Gamma)^{-1}\calB u_0](t)$ in the nonautonomous case.
By equation \eqref{nonauttransfcn},
this coincides with the frequency response function for time-varying
systems which is  defined in \cite[Corollary~3.2]{BGK}  by the formula
$$
\int_0^t C(t)U(t,\tau) B(\tau)u_0 e^{i\omega (\tau-t)}\,d\tau.
$$
Also, as noted in this reference, the result of our derivation agrees
 the Arveson frequency response function as
it appears in~\cite{SK}.  We recover it here explicitly as the Laplace
transform of an input-output operator (see equation \eqref{LapL}).

\subsection{Autonomous systems}\label{AUTsubsection} %-------------

In this subsection we give the proof of \eqref{rstabBanIneq}
 when  $X$, $U$ and $Y$ are Banach spaces.  In the process, however, we
also
consider two other ``stability radii": a pointwise stability radius and
a dichotomy radius.

First, we give a generalization to Banach spaces of
Theorem~\ref{GearTh} (cf.~\cite{LMS2}).  Here,
$\Fper$ denotes the Banach space $L^p([0,2\pi],X)$, $1\le
p<\infty$.
If $\{e^{tA}\}_{t\ge0}$ is a strongly continuous
semigroup on $X$, $\{E_{per}^t\}_{t\ge0}$
will denote the evolution semigroup defined on $\Fper$ by the rule
$E^t_{per}f(s)=e^{tA}f([s-t](\mod 2\pi))$;  its generator will be
denoted by $\Gammaper$.  The symbol
$\Lambda$ will be used to denote the set of all finite sequences
$\{v_k\}_{k=-N}^N$ in $X$ or $\calD(A)$, or $\{u_k\}_{k=-N}^N$ in $U$.

\begin{thm}\lb{LMS-Gearhart+BDeltaC}
Let $A$ generate a $C_0$ semigroup $\{e^{tA}\}_{t\ge0}$ on
$X$.  Let $B$ and $C$ be as above, and $\Delta\in \mathcal(Y,U)$.  Let
$\{e^{t(A+B\Delta C)}\}_{t\ge0}$ be the strongly continuous
semigroup  generated by $A+B\Delta C$.
Then the following are equivalent:
\begin{enumerate}
\item $1\in \rho (e^{2\pi (A+B\Delta C)})$;
\item $i\bbZ \subset \rho (A+B\Delta C)$ \quad and
  \[\disp{  \sup_{\{v_k\}\in\Lambda}
  \frac{\|\sum_k(A-ik+B\Delta C)^{-1}
  v_ke^{ik(\cdot)}\|_\Fper}{\|\sum_k v_k
  e^{ik(\cdot)}\|_\Fper} <\infty;}\]
\item $i\bbZ \subset \rho (A+B\Delta C)$ \quad and
  \[\disp{\inf_{\{v_k\}\in\Lambda}
  \frac{\|\sum_k(A-ik+B\Delta C)v_ke^{ik(\cdot)}\|_\Fper}{\|\sum_k v_k
  e^{ik(\cdot)}\|_\Fper} >0.}\]
\end{enumerate}
Further, if $\Gammaper$ denotes the generator of the evolution
 semigroup on $\Fper$, as above, and if $1\in\rho(e^{2\pi A})$, then
$\Gammaper$ is invertible and
\begin{equation}\lb{NormCGamInvB}
\|\calC\Gammaper^{-1}\calB\|=
\sup_{\{u_k\}\in\Lambda}
     \frac{\|\sum_k C(A-ik)^{-1}Bu_ke^{ik(\cdot)}\|_{L^p([0,2\pi],Y)}}
     {\|\sum_k u_ke^{ik(\cdot)}\|_{L^p([0,2\pi],U)}}
\end{equation}
where $\calC\Gammaper^{-1}\calB\in
\mathcal{L}(L^p([0,2\pi],U),L^p([0,2\pi],Y))$.
\end{thm}

\begin{proof}  The equivalence of {\em (i)--(iii)} follows as in
Theorem 2.3 of \cite{LMS2}.  For the last statement, let $\{u_k\}$ be a
finite set in $U$ and consider functions
$f$ and $g$ of the form
\[
f(s)=\sum_k(A-ik)^{-1}Bu_ke^{iks}\, ,\quad\mbox{and}\quad
g(s)=\sum_kBu_ke^{iks}\, .
\]
Then $f=\Gammaper^{-1}g$.  For,
\begin{align*}
(\Gammaper f)(s)&=
\left.\frac{d}{dt}\right|_{t=0}e^{tA}f([s-t]\mbox{mod}2\pi)\\
&= \sum_k [A(A-ik)^{-1}Bu_ke^{iks}-ik(A-ik)^{-1}Bu_ke^{iks}]=g(s).
\end{align*}
For functions of the form $h(s)=\sum_k u_ke^{iks}$, where
$\{u_k\}_k$ is a finite set in $U$, we have
$\calC\Gamma^{-1}\calB h=\sum_kC(A-ik)^{-1}Bu_ke^{ik(\cdot)}$.
Taking the supremum over all such functions gives:
\begin{align*}
\|\calC\Gammaper^{-1}\calB\|&=
\sup_{h}\frac{\|\calC\Gammaper^{-1}\calB h\|}{\|h\|}\\
&=\sup_{\{u_k\}\in\Lambda}
\frac{\|\sum_k C(A-ik)^{-1}Bu_ke^{ik(\cdot)}\|_{L^p([0,2\pi],Y)}}
{\|\sum_k u_ke^{ik(\cdot)}\|_{L^p([0,2\pi],U)}}.
\end{align*}
\end{proof}

In view of these facts we introduce a ``pointwise" variant of the
constant stability radius:
for $t_0>0$ and  $\lambda\in\rho(e^{t_0 A})$,
define the {\em pointwise stability radius}
\[
\rcstab^\lambda(e^{t_0 A},B,C):=
\sup\{r>0: \|\Delta\|_{\mathcal{L}(Y,U)}\le r
\Rightarrow \lambda\in\rho(e^{t_0(A+B\Delta C)})\}.
\]
By rescaling,
the study of this quantity can be reduced to the case of $\lambda=1$ and
$t_0=2\pi$.  Indeed,
\[\rcstab^\lambda(e^{t_0 A},B,C)=
\frac{2\pi}{t_0}\rcstab^\lambda(e^{2\pi A'},B,C),\quad
\mbox{ where }\quad A'=\frac{t_0}{2\pi}A.\]
Also, after writing
$\lambda=|\lambda|e^{i\theta}$
($\theta\in\bbR$), note that
\[\rcstab^\lambda(e^{2\pi A},B,C)=\rcstab^1(e^{2\pi A''},B,C),\quad
\mbox{ for }\quad
A''=A-\frac{1}{2\pi}(\ln|\lambda|+i\theta).\]
Therefore,
$$
\rcstab^\lambda(e^{t_0 A},B,C)=\frac{2\pi}{t_0}\rcstab^1(e^{2\pi
A'''},B,C)
$$
for
$$
A'''=\frac{1}{2\pi}(t_0A-\ln|\lambda|-i\theta).
$$

In the following theorem we estimate
$\rcstab^1(e^{2\pi A},B,C)$.  The idea for the proof goes back to
\cite{HP94}.
See also further developments in \cite{FishvN}.

\begin{thm}\lb{PtwiseBounds}
Let $\{e^{tA}\}_{t\ge0}$ be a strongly continuous semigroup generated
by $A$ on $X$, and assume $1\in\rho(e^{2\pi A})$.  Let $\Gammaper$
denote the
generator of the induced evolution semigroup on $\Fper$.
Let $B\in \mathcal{L}(U,X)$, and $C\in \mathcal{L}(X,Y)$. Then
\begin{equation}\lb{rstab1Bounds}
\frac{1}{\|\calC\Gammaper^{-1}\calB\|}
\le \rcstab^1(e^{2\pi A},B,C)
\le \frac{1}{\sup_{k\in\bbZ}\|C(A-ik)^{-1}B\|}.
\end{equation}
If $U$ and $Y$ are a Hilbert spaces and $p=2$, then equalities hold in
\eqref{rstab1Bounds}.
\end{thm}

\begin{proof}  The first inequality follows from an argument as in
Theorem~\ref{L0BoundThm}.  For the second inequality, let $\epsilon>0$,
and choose $\bar u\in U$ with $\|\bar u\|=1$ and $k_0\in\bbZ$ such that
\[
\|C(A-ik_0)^{-1}B\bar{u}\|_Y\ge
\sup_{k\in\bbZ}\|C(A-ik)^{-1}B\|-\epsilon>0.
\]
Using the Hahn-Banach Theorem, choose $y^*\in Y^*$ with
$\|y^*\|\le 1$ such that
\[
\left\langle y^*,\frac{C(A-ik_0)^{-1}B\bar{u}}
{\|C(A-ik_0)^{-1}B\bar{u}\|_Y}\right\rangle=1.
\]
Define $\Delta\in \mathcal{L}(Y,U)$ by
\[
\Delta y =
   -\frac{\langle y^*,y\rangle}{\|C(A-ik_0)^{-1}B\bar u\|_Y}\bar{u},
\quad y\in Y.
\]
We note that
\begin{equation}\lb{eqnA}
\Delta C(A-ik_0)^{-1}B\bar u=
-\frac{\langle y^*,C(A-ik_0)^{-1}B\bar u\rangle}
{\|C(A-ik_0)^{-1}B\bar{u}\|_Y}\bar u=-\bar{u},
\end{equation}
and
\begin{equation}\lb{eqnB}
\|\Delta\|\le\frac{1}{\|C(A-ik_0)^{-1}B\bar{u}\|_Y}
\le\frac{1}{\sup_{k\in\bbZ}\|C(A-ik_0)^{-1}B\bar{u}\|_Y-\epsilon}\, .
\end{equation}
Now set $\bar{v}:= (A-ik_0)^{-1}B\bar u$ in $X$.  By \eqref{eqnA},
$\Delta C\bar{v}=-\bar{u}$, and so
\[ (A-ik_0+B\Delta C)\bar{v}
=(A-ik_0)\bar{v}+B\Delta C\bar{v}=B\bar{u}- B\bar{u}=0.\]
Therefore,
\begin{equation*}
\inf_{\{v_k\}\in\Lambda}\frac{\|\sum_k(A-ik+B\Delta C)
   v_ke^{ik(\cdot)}\|_\Fper}{\|\sum_ku_ke^{ik(\cdot)}\|_\Fper}
\le \frac{\|(A-ik_0+B\Delta C)
   \bar{v}e^{ik_0(\cdot)}\|_\Fper}{\|\bar{v}e^{ik_0(\cdot)}\|_\Fper} =0.
\end{equation*}
By Theorem \ref{LMS-Gearhart+BDeltaC},
$1\notin\rho(e^{2\pi(A+B\Delta C)})$.  This shows that
$\rcstab^1(e^{2\pi A},B,C)\le\|\Delta\|$.

To finish the proof, suppose that
$ \rcstab^1(e^{2\pi A},B,C)>(\sup_{k\in\bbZ}\|C(A-ik)^{-1}B\|)^{-1}$.
Then with
$r:=(\sup_{k\in\bbZ}\|C(A-ik)^{-1}B\bar{u}\|_Y-\epsilon)^{-1}$,
 and $\epsilon>0$ chosen to be sufficiently small,
%%r:=\frac{1}{\sup_{k\in\bbZ}\|C(A-ik)^{-1}B\bar{u}\|_Y-\epsilon}\, ,
one has
\[ \frac{1}{\sup_{k\in\bbZ}\|C(A-ik)^{-1}B\bar{u}\|_Y}< r <
\rcstab^1(e^{2\pi A},B,C).
\]
But then by \eqref{eqnB}, $\|\Delta\|\le r < \rcstab^1(e^{2\pi A},B,C)$,
which is a contradiction.

For the last statement of the theorem, note that Parseval's formula
applied to
\eqref{NormCGamInvB} gives
\begin{equation}\label{ParForm}
\|\calC\Gammaper^{-1}\calB\|=\sup_{\{u_k\}\in\Lambda}
   \frac{\left(\sum_k\| C(A-ik)^{-1}Bu_k\|^2_Y\right)^{1/2}}
       {\left(\sum_k \|u_k\|^2_U\right)^{1/2}}
\le \sup_{k\in\bbZ}\|C(A-ik)^{-1}B\|.
\end{equation}
Therefore,
\[
\frac{1}{\|\calC\Gammaper^{-1}\calB\|}\ge
\frac{1}{\sup_{k\in\bbZ}\|C(A-ik)^{-1}B\|}
\]
and hence equalities hold in \eqref{rstab1Bounds}.
\end{proof}

Next we consider the following ``hyperbolic'' variant of the constant
stability radius. Recall, that a strongly continuous semigroup
$\{e^{tA}\}_{t\ge 0}$ on $X$ is called {\it hyperbolic} if
\[
\sigma(e^{tA})\cap\bbT=\emptyset,\quad\mbox{ where }
\bbT=\{z\in\bbC: |z|=1\},
\]
for some (and, hence, for all) $t>0$ (see, e.g., \cite{vanNbook}).
The hyperbolic semigroups are those for which the differential
equation $\dot x=Ax$ has exponential dichotomy (see, e.g., \cite{DK})
with the dichotomy projection $P$ being the Riesz projection
corresponding to the part of spectrum of $e^A$ that lies in
the open unit disc.
% with the dichotomy projection $P$ being the Riesz projection for
% the operator, say, $e^A$ that corresponds to the part of its
% spectrum contained in the open unit disc.

For a given hyperbolic semigroup $\{e^{tA}\}_{t\ge 0}$ and
 operators $B$, $C$ we define
 the {\it constant dichotomy radius} as:
\begin{equation*}
\begin{aligned}
rc_{dich}(\{e^{tA}\},B,C):= \sup\{
  r\ge 0&: \|\Delta\|_{\mathcal{L}(Y,U)}\le r \mbox{ implies }\\
  &\sigma(e^{t(A+B\Delta C)})\cap\bbT=\emptyset
  \mbox{ for all } t>0\}.
\end{aligned}
\end{equation*}
\noindent The dichotomy radius measures the size of the smallest
$\Delta\in\mathcal{L}(Y,U)$
for which  the perturbed equation $\dot x=[A+B\Delta C]x$ looses the
exponential dichotomy.

Now for any $\xi\in[0,1]$, consider the rescaled semigroup generated by
$A_\xi:=A-i\xi$ consisting of operators $e^{tA_\xi}=e^{-i\xi t}e^{tA}$,
$t\ge0$.  The pointwise stability radius can be related to the
dichotomy radius as follows.

\begin{lem}\label{dichrad}
 Let $\{e^{tA}\}_{t\ge0}$ be a hyperbolic semigroup.
Then
\[
rc_{dich}(\{e^{tA}\},B,C)=\inf_{\xi\in[0,1]}\rcstab^1(e^{2\pi
A_\xi},B,C).
\]
\end{lem}

\begin{proof} Denote the left-hand side by $\alpha$ and the right-hand
side by $\beta$.     First fix $r<\beta$.  Let $\xi\in[0,1]$. If
$\|\Delta\|\le r$, then $1\in\rho(e^{2\pi(A_\xi+B\Delta C)})$ and so
 $e^{i\xi 2\pi}\in\rho(e^{2\pi(A+B\Delta C)})$
for all $\xi\in[0,1]$.  That is, $e^{is}\in\rho(e^{2\pi(A+B\Delta C)})$
for all $s\in\bbR$, and so $\sigma(e^{2\pi(A+B\Delta C)})\cap\bbT
=\emptyset$.  This shows that $r\le\alpha$, and so
 $\beta\le\alpha$.

Now suppose $r<\alpha$.  If $\|\Delta\|\le r$, then
$\sigma(\{e^{t(A+B\Delta C)}\})\cap\bbT=\emptyset$, and so
$e^{i\xi t}\in\rho(e^{t(A+B\Delta C)})$ for all $\xi\in[0,1], t\in\bbR$.
That is, $1\in\rho(e^{t(A_\xi+B\Delta C)})$.
This says $r\le \beta$ and so $\alpha\le \beta$.
\end{proof}

Under the additional assumption that the semigroup
$\{e^{tA}\}_{t\ge 0}$ is exponentially stable
(that is, hyperbolic with a trivial dichotomy projection $P=I$),
Lemma~\ref{dichrad} gives, in fact, a formula for the constant {\it
stability} radius. Indeed, the following simple proposition holds.

\begin{prop}\label{dichst}
Let $\{e^{tA}\}_{t\ge0}$ be an exponentially stable semigroup.
Then
\[
rc_{dich}(\{e^{tA}\},B,C)= rc_{stab}(\{e^{tA}\},B,C).
\]
\end{prop}

\begin{proof} Denote the left-hand side by $\alpha$
and the right-hand side by $\beta$. Take $r<\beta$ and any $\Delta$
with $\|\Delta\|\le r$. By definition of the constant stability
radius, $\omega_0(\{e^{t(A+B\Delta C)}\})<0$. In particular,
$\sigma(e^{t(A+B\Delta C)})\cap\bbT=\emptyset$, and $r\le \alpha$ shows
that $\beta\le\alpha$.

Suppose that $\beta<r<\alpha$ for some $r$. By the definition of the
stability radius $\beta$,
there exists a $\Delta$ with $\|\Delta\|\in(\beta,r)$
such that
the semigroup $\{e^{t(A+B\Delta C)}\}_{t\ge0}$ is not stable.

For any $\tau\in[0,1]$
one has $\|\tau\Delta\|\le r<\alpha$. By the definition of the
dichotomy radius $\alpha$ it follows that the semigroup
$\{e^{t(A+\tau B\Delta C)}\}_{t\ge 0}$ is hyperbolic
for each $\tau\in[0,1]$.
Now consider its dichotomy projection
\[P(\tau)=(2\pi i)^{-1}\int_{\bbT}
\left(\lambda-e^{A+\tau B\Delta C}\right)^{-1}\,d\lambda,\]
which is the Riesz projection corresponding to the part
of $\sigma(e^{A+\tau B\Delta C})$ located inside of the open unit disk.
The function $\tau\mapsto P(\tau)$ is norm continuous.
Indeed, since the bounded perturbation $\tau B\Delta C$ of the
generator $A$ is continuous in $\tau$, the operators
$e^{t(A+\tau B\Delta C)}$, ${t\ge 0}$, depend on $\tau$
 continuously
(see, e.g., \cite[Corollary~3.1.3]{Pazy}); this implies
 the continuity of $P(\cdot)$ (see, e.g.,
\cite[Theorem~I.2.2]{DK}).

By assumption $\{e^{tA}\}_{t\ge 0}$ is exponentially stable, so
$P(0)=I$.
Also, $P(1)\neq I$ since the semigroup
$\{e^{t(A+B\Delta C)}\}_{t\ge 0}$ with
$\|\Delta\|\le r<\alpha$ is hyperbolic but not stable.
Since either  $\|I-P(\tau)\|=0$ or $\|I-P(\tau)\|\ge 1$, this
contradicts
the continuity of $\|P(\cdot)\|$.
\end{proof}

A review of the above development shows that the
inequality claimed in \eqref{rstabBanIneq} of Theorem~\ref{autin}
can now be proved.
\begin{proof}({\em of Theorem~\ref{autin}})
Indeed,
$\rstab(\{e^{tA}\},B,C)\le \rcstab(\{e^{tA}\},B,C)$, and so
\begin{alignat*}{3}
\frac{1}{\|\bbL\|}\le r_{stab}(\{e^{tA}\},B,C)
&\le  rc_{stab} (\{e^{tA}\},B,C)
       &\quad \mbox{(Theorem~\ref{L0BoundThm})} \\
\null & \le rc_{dich} (\{e^{tA}\},B,C)
       &\quad\mbox{(Proposition~\ref{dichst})} \\
\null & \le \inf_{\xi\in[0,1]}\rcstab^1(e^{2\pi A_\xi},B,C)
       &\quad \mbox{(Lemma~\ref{dichrad})} \\
\null &\le
\inf_{\xi\in[0,1]}\frac{1}{\sup_{k\in\bbZ}\|C(A_\xi-ik)^{-1}B\|}
       &\quad\mbox{(Theorem~\ref{PtwiseBounds})}\ \\
\null &= \frac{1}{\sup_{s\in\bbR}\|C(A-is)^{-1}B\|} &\null
\end{alignat*}
\end{proof}

We will need below the following simple corollary that holds for
{\it bounded} generators $A$.  (In fact, as shown in
\cite[Cor.~2.5]{FishvN}, formula \eqref{fBBA} below holds
provided $A$ generates a semigroup $\{e^{tA}\}_{t\ge0}$ that is
{\it uniformly} continuous just for $t>0$.)

\begin{cor}\label{bnddcase}
Assume $A\in{\mathcal L}(X)$ generates a (uniformly continuous) stable
semigroup on a Banach space $X$. Then
\begin{equation}\label{fBBA}
\rcstab(\{e^{tA}\},B,C)=\frac{1}{\sup_{s\in\bbR}\|C(A-is)^{-1}B\|}.
\end{equation}
\end{cor}
\begin{proof} By Theorem~\ref{autin}, it remains to prove only
the inequality ``$\ge$''. Fix $\Delta$ with $\|\Delta\|$ strictly less
than the right-hand side of \eqref{fBBA}.
Since $A+B\Delta C\in {\mathcal L}(X)$, it suffices
to show that $A+B\Delta C-\lambda=(A-\lambda)(I+(A-\lambda)^{-1}B\Delta
C)$
is invertible for each $\lambda$ with $\text{Re }\lambda\ge 0$. By the
analyticity of resolvent,
$\sup_{\text{Re }\lambda\ge
0}\|C(A-\lambda)^{-1}B\|\le\sup_{s\in{\mathbb
R}}\|C(A-is)^{-1}B\|$.
Thus,
$$
\begin{aligned}
\|\Delta\|&< \frac{1}{\sup_{s\in\bbR}\|C(A-is)^{-1}B\|}\\
&\le \frac{1}{\sup_{\text{Re }\lambda\ge 0}\|C(A-\lambda)^{-1}B\|}
\le \frac{1}{\|C(A-\lambda)^{-1}B\|},\quad \text{Re }\lambda\ge 0
\end{aligned}
$$
implies that $I+C(A-\lambda)^{-1}B\Delta$ is invertible. Therefore
(cf.~the proof of Theorem~\ref{L0BoundThm}), $I+(A-\lambda)^{-1}B\Delta
C$ is
invertible.
\end{proof}

\subsection{The norm of the input-output
operator}\label{normLsubsection}
%-----------

Since the lower bound on the stability radius is given by the norm of
the
input-output operator, which is defined by way of the solution
operators, it is of interest to express this quantity in terms of the
operators $A$, $B$ and $C$.  In this subsection it is shown that
for autonomous systems this
quantity can, in fact, be expressed explicitly in terms of the transfer
function:
\begin{equation}\lb{Lnorm}
\| \bbL \| = \sup_{u\in \calS(\bbR, U)}
\frac{\|\int_\bbR C(A-is)^{-1}Bu(s)e^{is(\cdot)}\,ds
\|_{L^p(\bbR,Y)}} {\|\int_\bbR u(s)e^{is(\cdot)}\,ds\|_{L^p(\bbR,U)}}.
\end{equation}
Here we use $\calS(\bbR, X)$ to denote the Schwartz class of
rapidly decreasing $X$-valued functions defined on $\bbR$: \
$\{v:{\mathbb R}\to X\big|\,\,\sup_{s\in{\mathbb
R}}\|s^mv^{(n)}(s)\|<\infty; n, m\in{\mathbb N}\}$.  As noted in
\eqref{rstabHilbeq}, $\|\bbL\|$ equals
$\sup_{s\in\bbR}\|C(A-is)^{-1}B\|$ if $U$ and $Y$ are Hilbert spaces and
$p=2$.
The section concludes by providing a similar expression, involving sums,
which serves as a lower bound for the constant stability radius.

The current focus is on autonomous systems so let $\{e^{tA}\}_{t\ge0}$
be a strongly continuous semigroup generated by $A$ and consider the
evolution semigroups $\{E^t_\bbR\}_{t\ge0}$ defined on functions on
the entire real line as
in \eqref{AUTevolsgR}, and $\{E^t\}_{t\ge0}$ defined for functions on
the
half-line as in \eqref{AUTevolsgHL}.  As before,  $\Gamma_\bbR$ and
$\Gamma$ will denote the generators of these semigroups on $L^p(\bbR,X)$
and $L^p(\bbR_+,X)$, respectively.  Both semigroups will be used as we
first show that $\|\calC\Gamma_\bbR^{-1}\calB\|$ equals the expression
in
\eqref{Lnorm} and then check that
$\|\bbL\|\equiv\|\calC\Gamma^{-1}\calB\|
=\|\calC\Gamma_\bbR^{-1}\calB\|$.

Given $v\in \calS(\bbR,X)$, let $g_v$  denote the function
$$
g_v(\tau)=\frac{1}{2\pi}\int_\bbR v(s) e^{i\tau s}\, ds,
\quad\tau\in{\mathbb R},
$$
and set $\scrG = \{g_v: v\in \calS(\bbR,X)\}.$
Assuming $\sup_{s\in\bbR}\|(A-is)^{-1}\|<\infty$,
define, for a given $v\in
\calS(\bbR,X)$, the function
$$f_v(\tau)=\frac{1}{2\pi}\int_\bbR (A-is)^{-1}v(s) e^{i\tau s}\, ds,
\quad\tau\in{\mathbb R},$$
and set $\scrF = \{f_v: v\in \calS(\bbR,X)\}$.

\begin{prop}\lb{FGprops}
Assume $\sup_{s\in\bbR}\|(A-is)^{-1}\|<\infty$. Then
\begin{enumerate}
\item $\scrG$ consists of differentiable functions, and is dense in
$L^p(\bbR,X)$;
\item $\scrF$ is dense in $\Dom(\Gamma_\bbR)$;
\item if $v\in \calS(\bbR,X)$ then $\Gamma_\bbR f_v= g_v$.
\end{enumerate}
\end{prop}

\begin{proof}
For $g\in L^1(\bbR, X)$, denote the Fourier transform by
\begin{equation*}
\hat{g} (\tau) = \frac{1}{2\pi}\int_\bbR e^{-is\tau} g(s)\, ds.
\end{equation*}
Note that $\scrG = \{g:\bbR \rightarrow X : \exists \ v\in \calS(\bbR,X)
\text{ so that } \hat{g}=v\}$, and so $\scrG$ contains the set $\{g\in
L^1(\bbR,X): \hat g\in \calS(\bbR,X) \}$. Since the latter set is dense
in $\LpR$, property (i) follows.

$\scrG$ consists of differentiable functions since for $ v\in
\calS(\bbR,X)$, the integral defining $g_v$ converges absolutely.
Moreover, for $ v\in \calS(\bbR,X)$, the function
$w(s)=(A-is)^{-1}v(s)$, $s\in \bbR$, is also in $\calS(\bbR,X)$,
since $\sup_{s\in{\mathbb R}}\|(A-is)^{-1}\|<\infty$. Hence $f_v$
is differentiable with derivative
\begin{equation*}
f_v'(\tau)=\frac{1}{2\pi}\int_\bbR is(A-is)^{-1}v(s) e^{i\tau s}\, ds=
\frac{1}{2\pi}\int_\bbR is w(s) e^{i\tau s}\, ds.
\end{equation*}
So $f'_v\in \LpR$, and hence $\scrF$ is dense in $\Dom(-d/dt + \calA)$.

Property (iii) follows from the following calculation:
\begin{align*}
(\Gamma f_v)(\tau) &= \frac{1}{2\pi}
\int_\bbR [-is(A-is)^{-1}v(s)e^{is\tau}+
A(A-is)^{-1}v(s)e^{is\tau}]\, ds\\
&= \frac{1}{2\pi}\int_\bbR (A-is)(A-is)^{-1}v(s)e^{is\tau}\, ds
= g_v(\tau).
\end{align*}
\end{proof}

Set $\Lambda_{\mathcal S}= \{v\in {\mathcal S}(\bbR,X): v(s)\in
\Dom (A) \mbox{ for }s\in \bbR,\ Av\in  {\mathcal S}(\bbR,X)  \}$.

\begin{prop}\lb{GammaRnorm}
Let $\{e^{tA}\}_{t\ge0}$ be a strongly continuous
semigroup generated by $A$.  Let
$\Gamma$ and $\Gamma_\bbR$ be the  generators of the
evolution semigroups on $\LpRp$ and $\LpR$, as defined in
\eqref{AUTevolsgHL} and \eqref{AUTevolsgR}, respectively.
 Then the following assertions hold:
\begin{enumerate}
\item if $\sigma(A)\cap i\bbR=\emptyset$ and
$\sup_{s\in{\mathbb R}}\|(A-is)^{-1}\|<\infty$
then $$\|\Gamma_\bbR\|_{\bullet,L^p(\bbR,X)}
=\inf_{v\in\Lambda_{\mathcal S}}
 \frac{\|\int_\bbR(A-is)v(s)e^{is(\cdot)}\,ds \| _\LpR}
 {\|\int_\bbR v(s)e^{is(\cdot)}\, ds\| _\LpR};$$
\item if $\Gamma_\bbR$ is invertible on $\LpR$, then
$\{e^{tA}\}_{t\ge0}$ is hyperbolic and
\begin{equation*}
\|\Gamma_\bbR^{-1}\|_{\mathcal{L}(\LpR)}=\sup_{v\in
\calS(\bbR,X)}\frac{\|\int_\bbR (A-is)^{-1}v(s)e^{is(\cdot)}\, ds \|
_\LpR}{\|
\int_\bbR v(s)e^{is(\cdot)}\, ds \| _\LpR};
\end{equation*}
\item if $\Gamma$ is invertible on $\LpRX$, then $\{e^{tA}\}_{t\ge 0}$
is
exponentially stable and $$\| \Gamma^{-1} \|_{\mathcal{L}(\LpRp)} =
\| \Gamma^{-1}_\bbR \|_{\mathcal{L}(\LpR)}.$$
\end{enumerate}
\end{prop}

\begin{proof}
To show (i) let $v\in \calS(\bbR,X)$. Since $\sup_{s\in{\mathbb
R}}\|(A-is)^{-1}\|<\infty$, the formula $w(s)=(A-is)^{-1}v(s),
 s\in \bbR$, defines a function, $w$, in $\Lambda_{\mathcal S}$.
Now,
\begin{equation*}
%\begin{aligned}
g_v(\tau) = \frac{1}{2\pi}
\int_\bbR (A-is)(A-is)^{-1}v(s)e^{is\tau}\,ds
 =\frac{1}{2\pi}\int_\bbR (A-is)w(s)e^{is\tau}\, ds
%\end{aligned}
\end{equation*}
and
\begin{equation*}
f_v(\tau)=\frac{1}{2\pi}\int_\bbR w(s)e^{is\tau}\, ds.
\end{equation*}
However, from Proposition \ref{FGprops},
\begin{equation*}
%\begin{align*}
\| \Gamma_\bbR \| _\bullet = \inf_{f_v\in \scrF}\frac{\|\Gamma_\bbR
f_v\|}{\|f_v \| }= \inf_{v\in \calS(\bbR,X)}\frac{\|g_v \|}{\|f_v\|}
=\inf_{w\in \Lambda_{\mathcal S}}
\frac{\| \int_\bbR (A-is)w(s)e^{is(\cdot)}\,
ds \| }{\| \int_\bbR w(s)e^{is(\cdot)}\, ds \| }\, .
%\end{align*}
\end{equation*}

To see (ii) note that
\begin{align*}
\| \Gamma_\bbR^{-1} \| &= \| \Gamma_\bbR \| _\bullet ^{-1} = \left[
\inf_{v\in \calS(\bbR,X)}\frac{\|\Gamma_\bbR f_v\|}{\|f_v \| }
\right]^{-1} =
%\left[ \inf_{v\in \calS}\frac{\|g_v \|}{\|f_v\|}  \right] ^{-1}\\ &=
\sup_{v\in \calS(\bbR,X)}\frac{\|f_v \|}{\|g_v\|}.
\end{align*}

For (iii) note that
$\| \Gamma_\bbR \| _{\bullet,\LpR} \le \| \Gamma\| _{\bullet, \LpRp}$.
  Indeed, let $f\in \LpR$ with $\supp f\subseteq \bbR_+$.  If
  $f\in \Dom(-d/dt + \mathcal{A})$, then $\supp \Gamma_\bbR f\subseteq
\bbR_+$ and $\|
\Gamma_\bbR f  \|_\LpR = \| \Gamma f \|_\LpRX$.
To see that $\| \Gamma_\bbR\| _\bullet \ge \| \Gamma \| _\bullet$\, ,
let $\epsilon>0$ and choose $f\in \Dom(-d/dt + \mathcal{A})$ with compact
support such that $\| f \|_\LpR =1$ and
$\| \Gamma_\bbR\| _\bullet \ge \|
\Gamma_\bbR f\| -\epsilon$.  Now choose $\tau \in \bbR$ such that
$f_\tau(s):= f(s-\tau)$, $s\in \bbR$, defines a function, $f_\tau \in
\LpR$, with $\supp f_\tau \subseteq \bbR_+$.  Let $\bar{f}_\tau$ denote
the element of $\LpRX$ which coincides with $f_\tau$ on $\bbR_+$.  Then
$\| f_\tau  \|=\| \bar{f}_\tau  \| $ and $\Gamma \bar{f}_\tau =-d/dt\,
f(\cdot -\tau) + Af(\cdot-\tau)=(\Gamma_\bbR f)_\tau$. Therefore, $\|
\Gamma_\bbR \| _\bullet \ge \| \Gamma_\bbR f \| -\epsilon = \|
(\Gamma_\bbR f )_\tau\| - \epsilon =  \| \Gamma_\bbR \bar{f} _\tau\|
-\epsilon \ge \| \Gamma \| _\bullet -\epsilon$.
\end{proof}

\begin{prop}\lb{GUdense}
The set $\scrG_U = \{ g_u: u\in \calS (\bbR,U)\}$ is dense in
$L^p(\bbR,U)$. If $u\in \calS (\bbR,U)$ and $B\in \mathcal{L}(U,X)$ then
$ B u\in
\calS (\bbR,X)$ and $\Gamma _\bbR f_{Bu}= \calB g_u$.
\end{prop}

\begin{proof}
The first statement is clear, as in Proposition \ref{FGprops}. The
second follows from the properties of Schwartz functions, and
from the calculation:
\begin{equation*}
\Gamma_\bbR f_{{\null}_{Bu}} = g_{{\null}_{Bu}}(\tau) =
\frac{1}{2\pi}\int_\bbR Bu(s) e^{is \tau}\, ds = B
\frac{1}{2\pi}\int_\bbR u(s) e^{is \tau}\, ds .
\end{equation*}
\end{proof}

Recall, see Remarks~\ref{AUTsmt} and Theorem~\ref{ExpStabThm+}, that
$\{e^{tA}\}_{t\ge 0}$ is hyperbolic (resp., stable) if and only if
$\Gamma_{\mathbb R}$
(resp., $\Gamma$) is invertible on $L^p({\mathbb R},X)$ (resp.,
$L^p({\mathbb
R}_+,X)$).

\begin{thm}\label{LnormThm}
If $\Gamma_\bbR$ is invertible on $\LpR$, then
\begin{equation}\lb{Lnorm1}
\| \calC  \Gamma_\bbR^{-1} \calB \| = \sup_{u\in \calS(\bbR, U)}
\frac{\|\int_\bbR C(A-is)^{-1}Bu(s)e^{is(\cdot)}\,ds
\|_{L^p(\bbR,Y)}} {\|\int_\bbR u(s)e^{is(\cdot)}\,ds\|_{L^p(\bbR,U)}}.
\end{equation}
If $\Gamma$ is invertible on $\LpRX$, then the norm of $\bbL = \calC
\Gamma^{-1} \calB $, as an operator from $\LpRU$ to $\LpRY$, is given by
the above formula:
\begin{equation}\label{Lnorm2}
\|\bbL\| = \| \calC \Gamma_\bbR^{-1} \calB \|.
\end{equation}
If, in addition, $U$ and $Y$ are Hilbert spaces and $p=2$, then
\begin{equation}\label{p=2}
\|\bbL\| = \sup_{s\in\bbR} \| C (A-is)^{-1} B
\|_{\mathcal{L}(U,Y)}.
\end{equation}
\end{thm}

\begin{proof}
For $u\in \calS(\bbR,U)$, consider functions $f_{Bu}$ and $g_u$.
Proposition \ref{GUdense} gives $f_{Bu}=\Gamma^{-1}_\bbR \calB g_u$ and
\begin{align*}
\| \calC \Gamma^{-1}_\bbR \calB \|
 &= \sup_{g_u\in \scrG_U}
\frac{\| \calC \Gamma^{-1}_\bbR \calB g_u \|_{L^p(\bbR,Y)}}
    {\| g_u  \|_{L^p(\bbR,U)}}
=\sup_{g_u\in \scrG_U}\frac{\|\calC f_{Bu}\|}{\| g_u  \|}\\
&=\sup_{u\in \calS(\bbR, U)}
\frac{\| \int_\bbR C(A-is)^{-1}Bu(s)e^{is(\cdot)}\,ds
\|_{L^p(\bbR,Y)}} {\|\int_\bbR u(s)e^{is(\cdot)}\,ds\| _{L^p(\bbR,U)}},
\end{align*}
which proves \eqref{Lnorm1}.

Now, if $\Gamma$ is invertible on $\LpRX$, then $\{e^{tA} \}_{t\ge 0}$
is
exponentially stable by Corollary~\ref{vanNCor}.
Hence, $\Gamma_\bbR$ is invertible on $\LpR$. Moreover,
for the case of the {\it stable} semigroup $\{e^{tA}\}_{t\ge 0}$,
the formula for $\Gamma_\bbR^{-1}$ (see, e.g., \cite{LaRa})
takes the form
\begin{equation*}
%\begin{align*}
(\Gamma_\bbR^{-1}f)(t)
=\int_0^\infty e^{sA}f(t-s)\, ds
%%  =  \int_{-\infty}^0 e^{-sA}f(t+s)\, ds
=\int_{-\infty}^t e^{(t-s)A}f(s)\, ds.
%% \quad f\in \LpR.
%\end{align*}
\end{equation*}

If $\supp f \subseteq (0,\infty)$, then
\begin{equation}\lb{GammaRinv}
(\Gamma_\bbR^{-1}f)(t)=\int_{-\infty}^t e^{(t-s)A}f(s)\, ds = \int_0^t
e^{(t-s)A}f(s)\, ds.
\end{equation}
For a function $h\in \LpRX$, define an extension $\tilde{h} \in \LpR$ by
$\tilde{h}(t) =
h(t)$ for $t\ge 0$ and $\tilde{h}(t)=0$ for $t<0$.
Then \eqref{GammaRinv} shows that $\Gamma_\bbR^{-1}\tilde{h}
=(\Gamma^{-1}h)\sptilde$. In particular, for $u\in \LpRU$,
$\widetilde{\bbL u}= \widetilde{\calC \Gamma^{-1}\calB u}=
\calC \Gamma_\bbR^{-1} \calB \tilde{u}$.
Therefore,
\begin{align*}
\| \bbL u \| _{\LpRY}
&= \| \widetilde{\bbL u}  \|_{L^p(\bbR,Y)}= \|
\calC \Gamma_\bbR^{-1}\calB \tilde{u}\|_{L^p(\bbR,Y)}\\
&\le \| \calC \Gamma_\bbR ^{-1} \calB \| \cdot
\|\tilde{u}\|_{L^p(\bbR,U)}=\| \calC \Gamma_\bbR ^{-1} \calB \|\cdot
\|u\|_{\LpRU}.
\end{align*}
This shows that $\| \bbL \| \le \| \calC \Gamma_\bbR ^{-1} \calB \|$.

To prove that equality holds in \eqref{Lnorm2}, let $\epsilon>0$
and choose $u\in L^p(\bbR,U)$, $\| u\| =1$, such that
$\|\calC \Gamma_\bbR ^{-1} \calB u \|_{L^p(\bbR,Y)}\ge
 \|\calC  \Gamma_\bbR^{-1}\calB \| -\epsilon$.
Without loss of generality, $u$ may be assumed to have compact support.
Now choose $r$ such that $\supp u(\cdot -r)\subseteq (0,\infty)$ and set
$w(\cdot):= u(\cdot -r)$.  Then $w\in L^p(\bbR,U)$ with $\supp
w\subseteq
(0,\infty)$.  Let $\bar{w}$ denote the element of $\LpRU$ that coincides
with $w$ on $\bbR_+$.  As in \eqref{GammaRinv} we have
\begin{equation*}
\calC \Gamma_\bbR ^{-1} \calB w(t)= C\int_0^t e^{(t-s)A} B w(s)\, ds= C
\int_{-\infty}^t e^{(t-s)A} B w(s)\, ds.
\end{equation*}
Since $\| \bar{w} \|_\LpRU = \| w \|_{L^p(\bbR,U)}=  \| u
\|_{L^p(\bbR,U)} =1$,  it follows that
\begin{align*}
\| \bbL \| &\ge \| \bbL \bar{w}  \| _\LpRY = \| \widetilde{\bbL \bar{w}}
\|_{L^p(\bbR,Y)}\\
&=\| \bbL \tilde{\bar{w}} \|_{L^p(\bbR,Y)}=
\| \calC \Gamma_\bbR ^{-1} \calB w\|_{L^p(\bbR,Y)}\\
%% &=\| C \int_{-\infty}^\cdot e^{(\cdot-s)A} B u(s-r)\, ds
%% \|_{L^p(\bbR,Y)}\\
%% &=\| C \int_{-\infty}^{\cdot -r}e^{(\cdot-r-\tau)A} B u(\tau)\, d\tau
%%  \|_{L^p(\bbR,Y)}\\
&=\| C \int_{-\infty}^\cdot e^{(\cdot-\tau)A} B u(\tau)\, d\tau \|
_{L^p(\bbR,Y)}\\
&=\| \calC \Gamma_\bbR ^{-1} \calB u \|_{L^p(\bbR,Y)}
\ge \| \calC \Gamma_\bbR ^{-1} \calB\| - \epsilon.
\end{align*}
This confirms \eqref{Lnorm2}. Parseval's formula and
\eqref{rstabBanIneq} give \eqref{p=2}.
\end{proof}

%edited to here 5-6-99

If $\{e^{tA}\}_{t\ge0}$ is exponentially stable
then the inequalities in \eqref{rstabBanIneq} give lower and upper
bounds on the stability radius in terms of $\bbL$ and $C(A-is)^{-1}B$,
respectively.  The previous theorem shows that $\|\bbL\|$ can be
explicitly expressed in terms of an integral involving $C(A-is)^{-1}B$.
We
conclude by observing that a lower bound for the constant stability
radius can be expressed by a similar formula involving a sum.
 For this, let $\xi\in[0,1]$ and
set
\[
S_\xi:=\sup_{\{u_k\}\in\Lambda}
     \frac{\|\sum_k
C(A-i\xi-ik)^{-1}Bu_ke^{ik(\cdot)}\|_{L^p([0,2\pi],Y)}}
     {\|\sum_k u_ke^{ik(\cdot)}\|_{L^p([0,2\pi],U)}}.
\]
We note that $S_\xi$ is computed as in equation \eqref{NormCGamInvB}
with $A$ replaced by $A_\xi=A-i\xi$.

\begin{cor}\lb{rstabBoundsCor}
Let $\{e^{tA}\}_{t\ge0}$ be an exponentially stable semigroup generated
by $A$.  Then
%\begin{equation}\lb{rstabBounds}
\[
\frac{1}{\sup_{\xi\in[0,1]}S_\xi}\le
\rcstab(\{e^{tA}\},B,C)
\le \frac{1}{\sup_{s\in\bbR}\|C(A-is)^{-1}B\|}.
\]
%\end{equation}
\end{cor}

\begin{proof}  Fix $\xi\in[0,1]$, and let $\Gammaperxi$ denote
the generator on $L^p([0,2\pi],X)$ of the evolution semigroup  induced
by
$\{e^{tA_\xi}\}_{t\ge0}$.  By Theorem \ref{LMS-Gearhart+BDeltaC},
$\|\calC\Gammaperxi^{-1}\calB\|=S_\xi$,
and so by Theorem \ref{PtwiseBounds},
\[
\frac{1}{S_\xi}\le
\rcstab^1(e^{2\pi A_\xi},B,C)
\le \frac{1}{\sup_{k\in\bbZ}\|C(A_\xi-ik)^{-1}B\|}.
\]
By Proposition~\ref{dichst}, taking the infimum over $\xi\in[0,1]$ gives
\begin{align*}
\frac{1}{\sup_{\xi\in[0,1]}S_\xi}&
\le \inf_{\xi\in[0,1]}\rcstab^1(e^{2\pi A_\xi},B,C)
=\rcstab(\{e^{tA}\},B,C)\\
&\le \inf_{\xi\in[0,1]}\frac{1}{\sup_{k\in\bbZ}\|C(A_\xi-ik)^{-1}B\|}\\
&=\frac{1}{\sup_{s\in\bbR}\|C(A-is)^{-1}B\|}.
\end{align*}
\end{proof}

\subsection{Two counterexamples}\label{counterex}

In contrast to the Hilbert space setting, the following Banach space
examples show that either inequality in \eqref{rstabBanIneq} may be
strict.
We start with the example where the second inequality in
\eqref{rstabBanIneq} is strict.

\begin{example}\label{RstabIneqEx} {\em
An example due to W.~Arendt (see, e.g., \cite{vanNbook},
Example~1.4.5) exhibits a (positive) strongly continuous
semigroup $\{e^{tA}\}_{t\ge 0}$ on a Banach space $X$
with the property that
$s_0(A)<\omega_0(A)<0$ for the abscissa of uniform boundedness of the
resolvent and the growth bound.  Now, for
$\alpha$ such that
$0\le\alpha\le-\omega_0(A)$, consider a rescaled semigroup  generated by
$A+\alpha$, and denote by $\Gamma_{A+\alpha}$ the generator of the
induced
evolution semigroup on $\LpRp$. The following relationships hold:
\[ \mbox{for } 0\le\alpha<-\omega_0(A),\quad
s_0(A+\alpha)=s_0(A)+\alpha<\omega_0(A)+\alpha=\omega_0(A+\alpha)<0;
\]
\[\mbox{for } \alpha_0:=-\omega_0(A),\quad
s_0(A+\alpha_0)<\omega_0(A+\alpha_0)=0.
\]
This says that  $s_0(A+\alpha)<0$ for all $\alpha\in[0,\alpha_0]$
and hence
\[
M:=\sup_{\alpha\in[0,\alpha_0]}
\sup_{s\in\bbR}\|(A+\alpha-is)^{-1}\|<\infty.
\]
Now note (see Corollary~\ref{vanNCor})
 that $\omega_0(A+\alpha)<0$ if and only if
$\|\Gamma_{A+\alpha}^{-1}\|<\infty$. Since $\omega_0(A+\alpha)\to 0$
as $\alpha\to\alpha_0$, we conclude that
$\|\Gamma_{A+\alpha}^{-1}\|\to\infty$ as $\alpha\to\alpha_0$.
Since  $\alpha\mapsto\|\Gamma_{A+\alpha}^{-1}\|$
is a continuous function of $\alpha$ on $[0,\alpha_0)$,
there exists $\alpha_1 \in [0,\alpha_0)$ such that
$\|\Gamma_{A+\alpha_1}^{-1}\|>M$, and so the following
inequality is strict:
\[
\frac{1}{\|\Gamma_{A+\alpha_1}^{-1}\|}
< \frac{1}{\sup_{s\in\bbR}\|(A+\alpha_1-is)^{-1}\|}.
\]
Also, we claim that there exists $\alpha_2\in[0,\alpha_0)$ such that
the following inequality is strict:
\[ \rcstab(\{e^{t(A+\alpha_2)}\},I,I) <
\frac{1}{\sup_{s\in\bbR}\|(A+\alpha_2-is)^{-1}\|}.
\]
To see this, let us suppose that for each $\alpha\in[0,\alpha_0)$ one
has $\rcstab(\{e^{t(A+\alpha)}\},I,I)\ge 1/(2M)$. Again, using that
$\omega_0(A+\alpha)\to 0$
as $\alpha\to\alpha_0$, find $\alpha\in[0,\alpha_0)$ such that
$|\omega_0(A+\alpha)|<1/(2M)$. Let $\Delta=\omega_0(A+\alpha)I$.
Since $\|\Delta\|=|\omega_0(A+\alpha)|$, by the definition
of stability radius one has:
\[0>\omega_0(A+\alpha+\Delta)=\omega_0(A+\alpha)-\omega_0(A+\alpha)=0,\]
a contradiction. Thus, there exists $\alpha_2\in[0,\alpha_0)$ such that
\[\rcstab(\{e^{t(A+\alpha_2)}\},I,I)\le \frac{1}{2M}<
\frac{1}{M}\le
\frac{1}{\sup_{s\in\bbR}\|(A+\alpha_2-is)^{-1}\|},\]
as claimed.\hfill$\diamondsuit$
}
\end{example}

This example shows that the second inequality in \eqref{rstabBanIneq}
can be strict due to the Banach-space pathologies related to the failure
of Gearhart's Theorem~\ref{GearTh}. Another example, given below,
shows that the first inequality in \eqref{rstabBanIneq} could be
strict due to the lack of Parseval's formula (see \eqref{ParForm}
in the proof of Theorem~\ref{PtwiseBounds}):
That is, the choice of $p=2$ in \eqref{rstabHilbeq} is as
important as the fact that $X$
in \eqref{rstabHilbeq} is a {\it Hilbert} space.
First, we need a formula for the norm of the input-output operator
on $L^1({\mathbb R}_+,X)$.

\begin{prop}\label{normonL1} Assume $\{e^{tA}\}_{t\ge0}$ is an
exponentially stable $C_0$ semigroup on a Banach space $X$.
The norm of the
operator ${\mathbb L}=\Gamma^{-1}$ on $L^1({\mathbb R}_+,X)$ is
\begin{equation}\label{normL1}
\|\Gamma^{-1}\|_{{\mathcal L}(L^1({\mathbb R}_+,X))}=
\sup_{\|x\|=1}\int\limits_0^\infty\|e^{tA}x\|\,dt.
\end{equation}
\end{prop}
\begin{proof}
Recall, see \eqref{AUTbbG}, that
\[\Gamma^{-1}f(t)=-\int\limits_0^t e^{\tau A}f(t-\tau)\,d\tau,
\quad t\in{\mathbb R}_+,\quad f\in L^1({\mathbb R}_+,X)\]
is the convolution operator. Choose positive
$\delta_n\in L^1({\mathbb R}_+,{\mathbb R})$ with
$\|\delta_n\|_{L^1}=1$ such that
\[
\|g*\delta_n-g\|_{L^1({\mathbb R}_+,X)}\to 0\quad\text{as}
\quad n\to\infty \quad\text{for each}\quad g\in L^1({\mathbb R}_+,X).
\]
Fix $x\in X$, $\|x\|=1$,
let $f=\delta_n x\in L^1({\mathbb R}_+,X)$ and note that
\[
\Gamma^{-1}f(t)=-\int\limits_0^te^{\tau A}x\,\delta_n(t-\tau)\,d\tau=
-(g*\delta_n)(t)\quad\text{for}\quad g(t)=e^{tA}x,\quad t\in{\mathbb
R}_+.
\]
This implies ``$\ge$'' in \eqref{normL1}. To see ``$\le$'', take
$f=\sum_{i=1}^N\alpha_i x_i$ with $\alpha_i\in L^1({\mathbb
R}_+,{\mathbb R})$ having disjoint supports and $\|x_i\|=1$,
$i=1,\ldots,N$.
Now $\|f\|_{L^1({\mathbb R}_+,X)}=\sum_i\|\alpha_i\|_{L^1}$ and for
$f_i(t)=e^{tA}x_i$ one has
\[\Gamma^{-1}f(t)=-\int\limits_0^t\sum_i e^{\tau
A}x_i\alpha_i(t-\tau)\,d\tau=-\sum_i(f_i*\alpha_i)(t).
\]
Using Young's inequality,
$$ \begin{aligned}
   \|\Gamma^{-1}f\|_{L^1({\mathbb R}_+,X)}
   &\le
   \sum_i\|f_i*\alpha_i\|_{L^1({\mathbb R}_+,X)}\\
   &\le
   \sum_i\|f_i\|_{L^1({\mathbb R}_+,X)}\|\alpha_i\|_{L^1} \le
\sup_{\|x\|=1}\int\limits_0^\infty\|e^{tA}x\|\,dt\sum_i\|\alpha_i\|_{L^1}.
\end{aligned}
$$
\end{proof}

\begin{example}\label{StephenEx}{\em
Take $X = {\mathbb C}^2$ with the $\ell_1$ norm.  Let
$$ A = \begin{pmatrix}
-1 & 1 \\ -1 & -1
\end{pmatrix}\quad
\text{so that}\quad
 e^{tA} = \begin{pmatrix}
e^{-t} \cos(t) & e^{-t} \sin(t) \\
-e^{-t} \sin(t) & e^{-t} \cos(t)
\end{pmatrix},$$
and
$$ (A-is)^{-1} = \frac{1}{(1+is)^2 + 1}
   \begin{pmatrix}
 -1-is & -1 \\ 1 & -1-is
\end{pmatrix} .$$
Since the extreme points of $X$ are $e^{i\theta} e_1$ and
$e^{i\theta} e_2$ ($\theta \in {\mathbb R}$), where
$e_1$ and $e_2$ are the unit vectors of ${\mathbb C}^2$, we see that
$$ \|(A-is)^{-1}\| = \frac{|1+is| + 1}{|(1+is)^2 + 1|} .$$
It may be numerically established that
$$ \sup_{s\in{\mathbb R}} \|(A-isI)^{-1}\| \approx 1.087494476 .$$
By Corollary~\ref{bnddcase}, the reciprocal to the last expression is
equal to $\rcstab(\{e^{tA}\},I,I)$.
On the other hand, using Proposition~\ref{normonL1},
$$ \begin{aligned}
   \|{\mathbb L}\| & =\|\Gamma_A^{-1}\|
   = \sup_{\|x\| = 1}
   \int_0^\infty \|e^{tA} x\| \, dt \\
   &= \int_0^\infty |e^{-t} \cos(t)| + |e^{-t} \sin(t)| \, dt
   \approx 1.262434309 .
\end{aligned}$$
Therefore, the first inequality in \eqref{rstabBanIneq} may be
strict.\hfill$\diamondsuit$  }\end{example}

The following example shows that the norm of the input-output
operator depends on $p$.
\begin{example}\label{StephenEx2}{\em
Let
$$ A = \begin{pmatrix}
9/2 & -5/2 \\ 25/2 & -13/2 \end{pmatrix} ,$$
acting on $\bbC^2$ with the Euclidean norm.  Thus
$$ e^{t A} = e^{-t} \begin{pmatrix}
   \cos t + (11/2)\sin t & -(5/2)\sin t \\
   (25/2) \sin t & \cos t - (11/2) \sin t \end{pmatrix} .$$
Then
$$ \|\Gamma^{-1}\|_{L_1\to L_1} \ge
   \int_0^\infty \|e^{tA} e_1\| \, dt
   \approx 7.748310791 , $$
whereas
$$ \|\Gamma^{-1}\|_{L_2\to L_2} =
   \sup_{s \in R} \|(A-is)^{-1}\|
   \approx 2.732492852 .$$}

\vspace{-19.2pt}\hfill$\diamondsuit$\vspace{19.2pt}
\end{example}

 %---------------------- section 4 -----------------------
\section{Internal and External Stability}

Work aimed at properties of stability and robustness of linear
time-invariant systems is often based on transform techniques.
More specifically, if the  transfer
function $H(\lambda)=C(A-\lambda)^{-1}B$ is a bounded analytic function
of
$\lambda$ in the right half-plane $\bbC_+=\{\lambda\in\bbC:\mbox{Re
}\lambda>0\}$, then the autonomous system
\eqref{mildAUTsys} is said to be {\em externally stable}.  This property
is often used to deduce {\em internal stability} of the system, i.e.,
the uniform exponential stability of the nominal system $\dot x=Ax$.
The relationship between internal and external stability has been
studied
extensively; see, e.g.,
\cite{And2, Cur, Cur2, JacNet, Log, Reb1, Reb2} and the references
therein.  In this section we examine the extent to which these
techniques apply to Banach-space settings and time-varying systems.  For
this, {\em input-output stability} of the system~\eqref{mildTVsys} will
refer to the property that the input-output operator $\bbL$ is  bounded
 from $L^p(\bbR_+,U)$ to $L^p(\bbR_+,Y)$.
If internal stability is assumed initially, then the
inequalities in \eqref{rstabBanIneq} exhibit a relationship between
these concepts of stability.
The next two theorems look at these relationships more closely and show,
in particular, when internal stability may be deduced from one of the
``external" stability conditions.  Therefore, throughout this section
$\{U(t,\tau)\}_{t\ge\tau}$ will denote a strongly continuous
exponentially bounded evolution family that is {\em not} assumed to be
exponentially stable.

\subsection{The nonautonomous case} In this subsection we give a very
short
proof of the fact that for general nonautonomous systems on Banach
spaces,
internal stability is equivalent to stabilizability, detectability and
input-output stability.  Before proceeding, it is worth reviewing some
properties of time-invariant systems.  For this, let
$\{e^{tA}\}_{t\ge0}$ be a
strongly continuous semigroup generated by $A$ on $X$, and let
$H_+^\infty(\mathcal{L}(X))$ denote the space of operator-valued
functions $G:\bbC\to\mathcal{L}(X)$ which are analytic on $\bbC_+$ and
$\sup_{\lambda\in\bbC_+}\|G(\lambda)\|<\infty$.    If $X$ is a
Hilbert space, it is well known that $\{e^{tA}\}_{t\ge0}$ is
exponentially stable if and only if $\lambda\mapsto (\lambda-A)^{-1}$ is
an
element of $H_+^\infty(\mathcal{L}(X))$; see, e.g., \cite{CZ},
Theorem~5.1.5.
This is a consequence of the fact that when $X$ is a Hilbert space,
$s_0(A)=\omega_0(e^{tA})$ (see \cite{vanNbook} or
Theorem~\ref{GearTh}).  If
$X$ is a Banach space, then strict inequality $s_0(A)<\omega_0(e^{tA})$
can
hold, and so exponential stability is no longer determined by the
operator $G(\lambda)=(\lambda-A)^{-1}$.  Extending these ideas to
address
systems~\eqref{mildAUTsys}, one considers
$H(\lambda)=C(\lambda-A)^{-1}B$:
it can be shown that if $U$ and $Y$ are Hilbert spaces, then {\em
\eqref{mildAUTsys} is internally stable if and only if it is
stabilizable, detectable and externally stable} (i.e.,
$H(\cdot)\in H_+^\infty(\mathcal{L}(U,Y))$).
See R.~Rebarber~\cite{Reb1} for a general result of this type.
It should be pointed out that
this work of Rebarber and others more recently
allows for a certain degree of
unboundedness of the operators $B$ and $C$.  Such ``regular" systems
(see
\cite{Weiss1}), and their time-varying generalizations, might be
addressed by combining the techniques of the present paper (including
the characterization of generation of evolution semigroups as found in
\cite{RRSV}) along and with those of \cite{HP94} and \cite{JDP}.  This
will not be done here.

If one allows for Banach
spaces, the conditions of stabilizability and detectability
are {\em not} sufficient to ensure that external stability  implies
internal stability.   Indeed, let $A$ generate a semigroup  for which
$s_0(A)<\omega_0(e^{tA})=0$ (see Example~\ref{RstabIneqEx}).  Then  the
system~\eqref{mildAUTsys} with
$B=I$ and $C=I$  is trivially stabilizable and detectable and
externally stable.  But since $\omega_0(e^{tA})=0$, it is not internally
stable.

Since the above italicized statement concerning external
stability fails for Banach-space systems~\eqref{mildAUTsys} and does not
apply to time-varying systems~\eqref{mildTVsys}, we aim to prove the
following extension of this.

\begin{thm}\label{IOstab} The system~\eqref{mildTVsys} is internally
stable if and only if it is stabilizable, detectable and input-output
stable.
\end{thm}

\noindent This theorem appears as part of Theorem~\ref{bigIOstab} below.
A version of it for finite-dimensional time-varying systems was
proven by B.~D.~O.~Anderson in \cite{And2}.
The fact that Theorem~\ref{IOstab} actually {\em extends} the
Hilbert-space statement above follows from the fact that the
 Banach-space inequality
$\sup_{\lambda\in \bbC_+}\|H(\lambda)\|\le \|\bbL\|$
(see~\cite{Weiss2}) which relates the operators that define external
and input-output stability is actually an {\em equality} for
Hilbert-space systems (see also~\cite{Weiss1}).

 In Theorems~\ref{IOstab} and \ref{bigIOstab}, below, the following
definitions are used.

\begin{defn}\label{StabDetDefn}
The nonautonomous system \eqref{mildTVsys}  is said to be
\begin{enumerate}
\item[(a)]  {\em stabilizable} if there
exists $F(\cdot)\in L^\infty(\bbR_+,\mathcal{L}_s(X,U))$ and a
corresponding exponentially stable evolution family
$\{U_{BF}(t,\tau)\}_{t\ge \tau}$ such that, for $t\ge s$ and $\ x\in
X$, one has:
\begin{equation}\lb{UBF}
U_{BF}(t,s)x=U(t,s)x+ \int_s^t
U(t,\t)B(\t)F(\t)U_{BF}(\t,s)x\,d\tau;
\end{equation}
\item[(b)] {\em detectable} if there
exists $K(\cdot)\in L^\infty(\bbR_+,\mathcal{L}_s(Y,X))$ and a
corresponding exponentially stable evolution family
$\{U_{KC}(t,\tau)\}_{t\ge \tau}$ such that, for $t\ge s$ and $\ x\in X$,
one has:
\begin{equation}\lb{UKC}
U_{KC}(t,s)x=U(t,s)x+ \int_s^t
U_{KC}(t,\t)K(\t)C(\t)U(\t,s)x\,d\tau .
\end{equation}
\end{enumerate}
An autonomous control system is called {\em stabilizable} if there is
an operator $F\in\mathcal{L}(X,U)$ such that $A+BF$ generates a
uniformly exponentially stable semigroup; that is, $\omega_0(A+BF)<0$.
Such a system is {\em detectable} if there is
an operator $K\in\mathcal{L}(Y,X)$ such that $A+KC$ generates a
uniformly exponentially stable semigroup.
\end{defn}

Using  Theorem~\ref{TVvanN}  to characterize exponential stability in
terms of the operator $\bbG$ as in \eqref{TVbbG} makes the proof of the
following theorem a straightforward manipulation of the appropriate
operators.

\begin{thm}\lb{bigIOstab}
The following are equivalent for a strongly continuous exponentially
bounded evolution family of operators $\U=\{U(t,\tau)\}_{t\ge \tau}$ on
a
Banach space $X$.

\newcounter{ii}
\begin{enumerate}
\item $\U$ is exponentially stable on $X$;
\item $\bbG$ is a bounded operator on $\LpRp$;
\item system \eqref{mildTVsys} is stabilizable and $\bbG\calB$ is a
bounded
operator from $\LpRU$ to $\LpRp$;
\item system \eqref{mildTVsys} is detectable and $\calC\bbG$ is a
bounded
operator from $\LpRp$ to $\LpRY$;
\item system \eqref{mildTVsys} is stabilizable and detectable and
$\bbL=\calC\bbG\calB$ is a bounded operator from $\LpRU$ to $\LpRY$.
\end{enumerate}
\end{thm}

\begin{proof}
The equivalence of (i) and (ii) is the equivalence of (i) and (ii) in
Theorem~\ref{TVvanN}.

To see that (ii) implies (iii), (iv), and (v), note that $\calB$ and
$\calC$
are  bounded, and thus $\bbL$ is bounded when $\bbG$ is bounded. So when
(ii) holds, the exponential stability of $\calU$, together with
boundedness
of $B(\cdot)$, $C(\cdot)$, $F(\cdot)$ and $K(\cdot)$, assure the
existence
of the evolution families $\{U_{BF}(t,\tau)\}_{t\ge \tau}$ and
$\{U_{KC}(t,\tau)\}_{t\ge \tau}$ as solutions of the integral equations
in
Definition ~\ref{StabDetDefn}; thereby showing that (iii), (iv), and (v)
hold.

To see that (iii) $\Rightarrow$ (ii), first note that the assumption of
stabilizability assures the existence of an exponentially stable
evolution
family $\mathcal{U}_{BF}=\{U_{BF}(t,\tau)\}_{t\ge \tau}$ satisfying
equation ~\eqref{UBF} for some $F(\cdot)\in
L^\infty(\bbR_+,\mathcal{L}_s(X,U))$. Given this exponentially stable
family, we define the operator $\bbG_{BF}$ by
\begin{equation}\lb{GBFdefn}
\bbG_{BF}f(s):= \int_0^s U_{BF}(s,\t)f(\t)\,d\t = \int_0^\infty
(E_{BF}^{\t}f)(s)\,d\t
\end{equation}
where $\{E_{BF}^{\t}f\}_{t\ge 0}$ is the semigroup induced by the
evolution
family $\calU_{BF}$ as described in equation ~\eqref{TVevolsgHL}.
$\bbG_{BF}$ is a bounded operator on $\LpRX$ by the equivalence of (i)
and
(ii).

For $f(\cdot) \in \LpRX $ and $s \in \bbR_+$, take $x=f(s)$ in equation
~\eqref{UBF}. Then, let $\xi = \t -s$, to obtain
\[
U_{BF}(t,s)f(s)=U(t,s)f(s)+ \int_0^{t-s}
U(t,\xi +s)B(\xi +s)F(\xi +s)U_{BF}(\xi +s,s)f(s)\,d\xi.
\] From this equation and from the definition of the semigroups $\{E^{t}\}_{t
\ge 0}$, and $\{E^{t}_{BF}\}_{t \ge 0}$ we obtain
\[
(E^{t-s}_{BF}f)(t) =(E^{t-s}f)(t) + \int_0^{t-s}(E^{t-s-\xi}\calB \calF
E_{BF}^{\xi}f)(t)\, d\xi
\]
and hence for $0 \le r$ and $0\le \sigma$ that
\[
(E^{r}_{BF}f)(\sigma) =(E^{r}f)(\sigma) + \int_0^{r}(E^{r-\xi}\calB
\calF
E_{BF}^{\xi}f)(\sigma)\, d\xi.
\]
Integrate from 0 to $\infty$ to obtain
\[
(\bbG_{BF }f)(\sigma) = (\bbG f)(\sigma) +
\int_0^{\infty}\int_0^{r}(E^{r-\xi}\calB \calF E_{BF}^{\xi}f)(\sigma)\,
d\xi\, dr.
\]
Let $r=\zeta +\eta$ and $\xi = \eta$ to obtain
\begin{equation}\lb{GBF}
\begin{aligned}
(\bbG_{BF }f)(\sigma) &= (\bbG f)(\sigma) +
\int_0^{\infty}\int_0^{\infty}(E^{\zeta}\calB \calF
E_{BF}^{\eta}f)(\sigma)\, d\eta\, d\zeta\\
&= (\bbG f)(\sigma) + (\bbG \calB \calF \bbG_{BF} f)(\sigma).
\end{aligned}
\end{equation}
That $\bbG$ is bounded now follows from equation ~\eqref{GBF}, the
boundedness of $\bbG\calB$, and the boundedness of $\bbG_{BF}$ and
$\calF$.

To see that (iv) $\Rightarrow$ (ii), first note that the assumption of
detectability assures the existence of an exponentially stable evolution
family $\mathcal{U}_{KC}=\{U_{KC}(t,\tau)\}_{t\ge \tau}$ satisfying
equation ~\eqref{UKC} for some $K(\cdot)\in
L^\infty(\bbR_+,\mathcal{L}_s(Y,X))$. Given this exponentially stable
family, the operator $\bbG_{KC}$, defined in a manner analogous to
$\bbG_{BF}$ in equation ~\eqref{GBFdefn}, is a bounded operator on
$\LpRX$.
A derivation beginning with equation ~\eqref{UKC}, and similar to that
which gave equation ~\eqref{GBF}, now gives $\bbG_{KC}  = \bbG  +
\bbG_{KC}\calK\calC\bbG$.  This equation, together with the assumed
boundedness of $\bbG_{KC}$, $\calK$, and $\calC \bbG$, gives the
boundedness of $\bbG$.

Finally, to see that (v) $\Rightarrow$ (ii), again note that the
assumption
of detectability yields an exponentially stable evolution family
$\calU_{KC}$ and an associated bounded operator $\bbG_{KC}$. For
$u(\cdot)
\in \LpRU$, and $s\in \bbR_+$ take $x = B(s)u(s)$ in equation
~\eqref{UKC}.
A calculation similar to that which gave equation ~\eqref{GBF} now gives
$\bbG_{KC}\calB  = \bbG\calB  + \bbG_{KC}\calK\calC\bbG\calB$.  The
assumed
boundedness of $\bbL=\calC \bbG \calB$, $\calK$, and $\bbG_{KC}$, now
yields the boundedness of $\bbG \calB$.  The boundedness of $\bbG \calB$
together with the assumption of stabilizability implies that $\bbG$ is
bounded by the equivalence of  (iii) and (ii).
\end{proof}

\subsection{The autonomous case}%----------------------------------

The main result of this subsection
is Theorem~\ref{autonstab} which builds on Theorem~\ref{LnormThm}
and parallels Theorem~\ref{bigIOstab}
for autonomous systems of the form \eqref{mildAUTsys}.
The main point is to provide explicit conditions, in terms of the
operators $A$, $B$ and $C$, which imply internal stability.

Let $A_\alpha := A-\alpha I$ denote the
generator of the rescaled semigroup
$\{e^{-\alpha t} e^{tA}\}_{t\ge 0}$.

\begin{thm}\lb{autonstab}
Let $\{e^{tA}\}_{t\ge 0}$ be a strongly continuous semigroup on
a Banach space $X$ generated by $A$. Let $U$ and $Y$ be Banach spaces
and assume $B\in\mathcal{L}(U,X)$ and $C\in\mathcal{L}(X,Y)$.  Then
the following are equivalent.
% \newcounter{numb}
\begin{enumerate}
% {\upshape(\arabic{numb})}{\usecounter{numb}\setlength{\labelsep}{10pt}
% \setlength{\itemsep}{10pt}\setlength{\leftmargin}{24pt}}
\item $\{e^{tA}\}_{t\ge 0} $ is exponentially stable;
\newline
\item $\bbG$ is a bounded operator on $\LpRX$;
\newline
\item $\sigma(A)\cap\overline{\mathbb C}_+=\emptyset$ and
$\displaystyle\sup_{v\in \calS(\bbR,X)}
\frac{\|\int_\bbR
  (A_\alpha-is)^{-1}v(s)e^{is(\cdot)}\,ds\|_{L^p(\bbR,X)}}
  {\|\int_\bbR v(s)e^{is(\cdot)}\,ds\|_{L^p(\bbR,X)}}<\infty$ \\for all
  $\alpha \ge 0$;
\newline
\item $\sigma(A)\cap\overline{\mathbb
C}_+=\emptyset$, $\displaystyle\sup_{u\in \calS(\bbR,U)}
\frac{\|\int_\bbR
  (A_\alpha-is)^{-1}Bu(s)e^{is(\cdot)}\,ds\|_{L^p(\bbR,X)}}
{\|\int_\bbR u(s)e^{is(\cdot)}\,ds\|_{L^p(\bbR,U)}}<\infty$
\\[5pt]
for all $\alpha \ge 0$, and  \eqref{mildAUTsys} is stabilizable;
\newline
\item $\sigma(A)\cap\overline{\mathbb C}_+=\emptyset$,
$\displaystyle\sup_{v\in \calS(\bbR,X)}
\frac{\|\int_\bbR
C(A_\alpha-is)^{-1}v(s)e^{is(\cdot)}\,ds\|_{L^p(\bbR,Y)}}
{\|\int_\bbR v(s)e^{is(\cdot)}\,ds\|_{L^p(\bbR,X)}}<\infty$
\\[5pt]
for all $\alpha \ge 0$,   and \eqref{mildAUTsys} is detectable;
\newline
\item $\sigma(A)\cap\overline{\mathbb
C}_+=\emptyset$, $\displaystyle\sup_{u\in \calS(\bbR,U)}
\frac{\|\int_\bbR
  C(A_\alpha-is)^{-1}Bu(s)e^{is(\cdot)}\,ds\|_{L^p(\bbR,Y)}}
{\|\int_\bbR u(s)e^{is(\cdot)}\,ds\|_{L^p(\bbR,U)}}<\infty$ \\[5pt]
for all $\alpha \ge 0$, and \eqref{mildAUTsys} is both stabilizable and
detectable.
\end{enumerate}
Moreover, if $\{e^{tA}\}_{t\ge0}$ is exponentially stable,
then the norm of the input-output operator,
$\bbL= \calC\bbG \calB $, is equal to
\[\sup_{u\in \calS(\bbR,U)}
\frac{\|\int_\bbR C (A-is)^{-1}B
u(s)e^{is(\cdot)}\,ds\|_{L^p(\bbR,Y)}}
{\|\int_\bbR u(s)e^{is(\cdot)}\,ds\|_{L^p(\bbR,U)}}.\]
\end{thm}
\begin{proof} First note the equivalence of statements
 (i) and (ii) follows from Remarks~\ref{AUTvanN}.
Also,  the implication (i)$\Rightarrow$(vi), as well as the last
statement
of the theorem, follow from  Theorem~\ref{LnormThm}.
The exponential stability of $\{ e^{tA}\}_{t\ge0}$ is equivalent
to the invertibility of $\Gamma_\bbR$ (Corollary~\ref{vanNCor}), and so
(iii) follows from (i) by Proposition~\ref{GammaRnorm}.

To see that (iii) implies (i), begin by setting $\alpha=0$.
We wish to use  properties of $\scrF$ as in
Proposition~\ref{FGprops}.  We begin by observing that if the
expression in
(iii) is finite, then
$\sup_{s\in\bbR}\|(A-is)^{-1}\|<\infty$.  Indeed, if this were not the
case, then there would exist $s_n\in\bbR$ and $x_n\in \text{Dom}(A)$
with
$\|x_n\|=1$ such that $\|(A-is_n)x_n\|\to 0$ as $n\to\infty$.
Choose functions $\beta_n\in\calS(\bbR)$ with the property that
\begin{equation}\label{betaneq}
\lim_{n\to\infty}\frac{\|\int_\bbR
  \beta_n(s)(is_n-is)e^{is(\cdot)}\,ds\|_{L^p(\bbR)}}
  {\|\int_\bbR \beta_n(s)e^{is(\cdot)}\,ds\|_{L^p(\bbR)}}=0.
\end{equation}
Note: to construct such a sequence of functions $\beta_n$,
one takes, without loss of generality, $s_n=0$ in \eqref{betaneq}
and chooses a ``bump" function $\beta_0(s)$ where $\beta_0(0)=1$ and
$\beta_0$ has support in $(-1,1)$.  Then set $\beta_n(s)=n\beta_0(ns)$.
If  $\check{\ }$ denotes the inverse Fourier transform, then
$\check{\beta}_n(\tau)=\check{\beta}_0(\tau/n)$.  Also, for
$\alpha_n(s)=s\beta_n(s)$, one has $\check{\alpha}_n(\tau)
= \frac{1}{n}\check{\alpha}_0(\tau/n)$, and so
$$
\frac{\|\int_\bbR\beta_n(s)se^{is\cdot}\,ds\|^p_{L^p(\bbR)}}
{\|\int_\bbR\beta_n(s)e^{is\cdot}\,ds\|^p_{L^p(\bbR)}}
=\frac{\|\check{\alpha}_n\|^p}{\|\check{\beta}_n\|^p}
=\frac{\left(\frac{1}{n}\right)^p
\|\check{\alpha}_0\|^p}{\|\check{\beta}_0\|^p}
\to 0,\quad\mbox{as }n\to\infty.
$$

Now, setting $v_n(s):=\beta_n(s)(A-is)x_n$ gives a function $v_n$ in
$\calS(\bbR,X)$ with the properties that
$(A-is)^{-1}v_n(s)=\beta_n(s)x_n$.
Thus,
\begin{align*}
&\frac{\|\int_\bbR(A-is)^{-1}v_n(s)e^{is(\cdot)}\,ds\|}
  {\|\int_\bbR v_n(s)e^{is(\cdot)}\,ds\|}
  =\frac{\|\int_\bbR
  \beta_n(s)x_ne^{is(\cdot)}\,ds\|}
  {\|\int_\bbR \beta_n(s)(A-is)x_n e^{is(\cdot)}\,ds\|}\\
&=\frac{\|\int_\bbR
  \beta_n(s)x_ne^{is(\cdot)}\,ds\|}
  {\|\int_\bbR \beta_n(s)(A-is_n)x_ne^{is(\cdot)}\,ds
    +\beta_n(s)(is_n-is)x_n e^{is(\cdot)}\,ds \|}\\
&\ge\frac{\|\int_\bbR
  \beta_n(s)x_ne^{is(\cdot)}\,ds\|}
  {\|(A-is_n)x_n\|\,\|\int_\bbR \beta_n(s)e^{is(\cdot)}\,ds\|
    +\|\int_\bbR\beta_n(s)(is_n-is)x_n e^{is(\cdot)}\,ds \|}\\
&=\left(\|(A-is_n)x_n\|+\frac{\|\int_\bbR
  \beta_n(s)(is_n-is)e^{is(\cdot)}\,ds\|}
  {\|\int_\bbR \beta_n(s)e^{is(\cdot)}\,ds\|}\right)^{-1}.
\end{align*}
By the choice of $s_n$, $x_n$
and $\beta_n$, this last expression goes to $\infty$ as $n\to\infty$,
contradicting (iii).  Hence if the expression in (iii) is finite for
$\alpha=0$, then
$\sup_{s\in\bbR}\|(A-is)^{-1}\|<\infty$.

Now we may apply
Proposition~\ref{FGprops} (ii) and Proposition~\ref{GammaRnorm}, to
obtain
$$
%\begin{align*}
\| \Gamma_\bbR \| _\bullet
= \inf_{v\in \calS (\bbR,X)}\frac{\|\Gamma_\bbR f_v\|} {\|f_v \| }=
\inf_{v\in \calS(\bbR,X)}\frac{\|g_v \|}{\|f_v\|}
= \left(\sup_{v\in \calS(\bbR,X)}\frac{\|f_v \|}{\|g_v\|}\right) ^{-1}
>0.
%\end{align*}
$$
This shows that $0\notin \sigma_{ap}(\Gamma_\bbR)$ and so, by
\cite{LMS2}, it follows that $\sigma_{ap}(e^{tA})\cap \bbT = \emptyset$.
On the other hand, since $\sigma(A) \cap i\bbR =\emptyset$, it follows
from the spectral mapping theorem for the residual spectrum,
$\sigma_r(e^{tA})$, that
\begin{equation*}
\sigma(e^{tA})\cap \bbT= \left[ \sigma_{ap}(e^{tA}) \cup
\sigma_{r}(e^{tA}) \right] \cap \bbT = \emptyset.
\end{equation*}
The same argument holds for any $\alpha \ge 0$. As a result,
$\{e^{tA_\alpha} \}_{t\ge0}$ is hyperbolic for each $\alpha\ge 0$,
and thus
$\{e^{tA}\}_{t\ge0}$ is exponentially stable.

So far it has been shown that the statements (i)--(iii) are
equivalent, and that statement (i) implies (vi).  By showing that
(vi)$\Rightarrow$(iv)$\Rightarrow$(iii)
and (vi)$\Rightarrow$(v)$\Rightarrow$(iii), we complete the proof.

To see that (vi) implies (iv),  begin by setting
$\alpha = 0$.  Since \eqref{mildAUTsys} is detectable, there exists
$K\in
\mathcal{L}(Y,X)$ such that $A + KC$ generates an exponentially stable
semigroup. By the implication  (i)$\Rightarrow$(iii) for the semigroup
$\{e^{t(A+KC)}\}$, it follows that
\begin{equation*}
M_1 := \sup_{v\in \calS(\bbR, X)}
\frac{\| \int_\bbR (A+KC-is)^{-1}v(s)e^{is(\cdot)}\,ds\| }
{\| \int_\bbR v(s)e^{is(\cdot)}\,ds\| }
\end{equation*}
is finite.  So,
\begin{equation}\lb{ineq1}
\begin{align}
&\sup_{u\in \calS(\bbR, U)}
\frac{\|\int_\bbR (A+KC-is)^{-1}Bu(s)e^{is(\cdot)}\,ds\|}
{\|\int_\bbR u(s)e^{is(\cdot)}\,ds\|}\\
&=\sup_{u\in \calS(\bbR, U)}
\frac{\|\int_\bbR (A+KC-is)^{-1}Bu(s)e^{is(\cdot)}\,ds\|}{\|\int_\bbR
Bu(s)e^{is(\cdot)}\,ds\|}\cdot\frac{\| B\int_\bbR
u(s)e^{is(\cdot)}\,ds\|}{\|\int_\bbR u(s)e^{is(\cdot)}\,ds\|}
\le M_1 \| B \|.\nonumber
\end{align}
\end{equation}
By hypothesis in (vi),
\begin{equation*}
M_2:=\sup_{u\in \calS(\bbR, U)}
\frac{\|\int_\bbR C(A-is)^{-1}Bu(s)e^{is(\cdot)}\,ds\|}
{\|\int_\bbR u(s)e^{is(\cdot)}\,ds\|}
\end{equation*}
is finite.  For $u\in \calS(\bbR,U)$,
let $w(s)=KC(A-is)^{-1}Bu(s),\quad s\in \bbR$. Then,
\begin{equation}\lb{ineq2}
\frac{\|\int_\bbR (A+KC-is)^{-1}KC(A-is)^{-1}Bu(s)e^{is(\cdot)}\,ds\|}
{\|\int_\bbR u(s)e^{is(\cdot)}\,ds\|}
\end{equation}
\begin{align*}
&=\frac{\|\int_\bbR (A+KC-is)^{-1}w(s)e^{is(\cdot)}\,ds\|}
{\|\int_\bbR w(s)e^{is(\cdot)}\,ds\|}
 \cdot \frac{\| K\int_\bbR C(A-is)^{-1}Bu(s)e^{is(\cdot)}\,ds\|}
{\|\int_\bbR u(s)e^{is(\cdot)}\,ds\|}\\
&\le M_1 \|K \| M_2.
\end{align*}

Finally, since
\begin{equation*}\begin{aligned}
(A-is)^{-1}B & = (A+KC-is)^{-1}B \\
& \qquad\qquad+ (A+KC-is)^{-1}KC(A-is)^{-1}B,
\end{aligned}
\end{equation*}
it follows from \eqref{ineq1} and \eqref{ineq2} that
\begin{equation*}
\sup_{u\in \calS(\bbR, U)}
\frac{\|\int_\bbR (A-is)^{-1}Bu(s)e^{is(\cdot)}\,ds\|}
{\|\int_\bbR u(s)e^{is(\cdot)}\,ds\|}\le M_1\| B \| + M_1 \| K \| M_2 .
\end{equation*}
This argument holds for all $\alpha\ge0$, so the implication
(vi)$\Rightarrow$(iv) follows.

To see that (iv) implies (iii) we again argue only in the case
$\alpha=0$. Since \eqref{mildAUTsys} is stabilizable, there exists $F\in
B(X,U)$ such that $A+BF$ generates an exponentially stable semigroup.
By the implication (i)$\Rightarrow$(iii) for the semigroup
$\{e^{t(A+BF)}\}$,
 it follows that
\begin{equation*}
M_3 := \sup_{v\in \calS(\bbR, X)}
\frac{\| \int_\bbR (A+BF-is)^{-1}v(s)e^{is(\cdot)}\,ds\| }
{\| \int_\bbR v(s)e^{is(\cdot)}\,ds\| }
\end{equation*}
is finite.  By hypotheses in (iv),
\begin{equation*}
M_4:=\sup_{u\in \calS(\bbR, U)}
\frac{\|\int_\bbR (A-is)^{-1}Bu(s)e^{is(\cdot)}\,ds\|}
{\|\int_\bbR u(s)e^{is(\cdot)}\,ds\|}
\end{equation*}
is finite.  For $v\in \calS(\bbR,X)$,
set $w(s)= F(A+BF-is)^{-1}v(s)$,
$s\in\bbR$. Then,
\begin{equation}\lb{ineq3}
\sup_{v\in \calS(\bbR, X)}
\frac{\| \int_\bbR(A-is)^{-1}BF (A+BF-is)^{-1}v(s)e^{is(\cdot)}\,ds\| }
{\| \int_\bbR v(s)e^{is(\cdot)}\,ds\| }
\end{equation}
\begin{align*}
&=\frac{\| \int_\bbR(A-is)^{-1}Bw(s)e^{is(\cdot)}\,ds\| }
{\| \int_\bbR w(s)e^{is(\cdot)}\,ds\| }\cdot\frac{\|
F\int_\bbR(A+BF-is)^{-1}v(s)e^{is(\cdot)}\,ds\| }
{\| \int_\bbR v(s)e^{is(\cdot)}\,ds\| }\\
&\le M_4 \| F \| M_3.
\end{align*}
Since
\begin{equation*}
(A-is)^{-1} = (A+BF-is)^{-1} + (A-is)^{-1}BF(A+BF-is)^{-1},
\end{equation*}
it follows from \eqref{ineq3} that
\begin{equation*}
\sup_{v\in \calS(\bbR, X)}\frac{\|
\int_\bbR(A-is)^{-1}v(s)e^{is(\cdot)}\,ds\| }
{\| \int_\bbR v(s)e^{is(\cdot)}\,ds\| }\le M_3 + M_4 \| F \| M_3.
\end{equation*}
Thus (iii) follows from (iv).
Similar arguments show that (vi)$\Rightarrow$(v) and
(v)$\Rightarrow$(iii).
\end{proof}
{}From the equivalence of statements (i) and (iii), it follows that
the growth bound of a semigroup on a Banach space is given by
$$
\omega_0(e^{tA})=\inf\left\{\alpha\in\bbR:\sup_{v\in \calS(\bbR,X)}
\frac{\|\int_\bbR
  (A_\alpha-is)^{-1}v(s)e^{is(\cdot)}\,ds\|}
  {\|\int_\bbR v(s)e^{is(\cdot)}\,ds\|}<\infty\right\}.
$$
This is a natural generalization of the formula for the growth bound
for a semigroup on a Hilbert space as provided by Gearhart's
Theorem, see \cite{Huang,Nagel,vanNbook,Prus} and cf.
Theorem~\ref{GearTh}:
$$
\omega_0(e^{tA})=s_0(A)=\inf\left\{\alpha\in\bbR:
\sup_{\mbox{Re}\lambda\ge\alpha}\|(A-\lambda)^{-1}\|<\infty\right\}.
$$

\begin{thm}\label{Rebarbthm}
Let $\{e^{tA}\}_{t\ge0}$ be a strongly continuous
 semigroup on a Banach space $X$ with the property that
$s_0(A)=\omega_0(e^{tA})$.  Assume \eqref{mildAUTsys} is stabilizable
and
detectable.
If $\overline\bbC_+\subset\rho(A)$ and
$M:=\sup_{s\in\bbR}\|C(A-is)^{-1}B\|<\infty$, then
$\{e^{tA}\}_{t\ge0}$ is exponentially stable.
\end{thm}

\begin{proof}
Choose operators $F\in \mathcal{L}(X,U)$ and $K\in \mathcal{L}(Y,X)$
such that the semigroups generated by $A+BF$ and $A+KC$ are
exponentially stable.  Then  $s_0(A+BF)<0$
and $s_0(A+KC)<0$, and so
$$
M_1:=\sup_{s\in\bbR}\|(A+BF-is)^{-1}\|,\mbox{ \ and \ }
M_2:=\sup_{s\in\bbR}\|(A+KC-is)^{-1}\|.
$$
are both finite.  Since
$$
(A-is)^{-1}B=(A+KC-is)^{-1}B+(A+KC-is)^{-1}KC(A-is)^{-1}B,
$$
it follows that
$$
M_3:=\sup_{s\in\bbR}\|(A-is)^{-1}B\|\le M_2\|B\|+M_2\|K\|M.
$$
Also,
$$
(A-is)^{-1}=(A+BF-is)^{-1}+(A-is)^{-1}BF(A+BF-is)^{-1},
$$
and so
$$
\sup_{s\in\bbR}\|(A-is)^{-1}\|\le M_1+M_3\|F\|M_1.
$$
Therefore, $\omega_0(e^{tA})=s_0(A)<0$.
\end{proof}

The following result, based on~\cite{KaLu},
describes a particular situation in
which $s_0(A)=\omega_0(e^{tA})$.

\begin{cor}\label{corKaLun} Assume that for the
generator $A$ of a strongly continuous
semigroup $\{e^{tA}\}_{t\ge0}$ on a Banach space $X$
there exists an $\omega>\omega_0(e^{tA})$ such that
\begin{equation}\label{condKaLu1}
\int\limits_{-\infty}^\infty \|
(\omega+i\tau-A)^{-1}x\|_X^2\,d\tau<\infty
\quad\text{for all}\quad x\in X,
\end{equation}
and
\begin{equation}\label{condKaLu2}
\int\limits_{-\infty}^\infty \|
(\omega+i\tau-A^*)^{-1}x^*\|_{X^*}^2\,d\tau<\infty
\quad\text{for all}\quad x^*\in X^*,
\end{equation}
where $X^*$ is the adjoint space. Then system
\eqref{mildAUTsys} is internally stable if and
only if it is  stabilizable, detectable and
externally stable.
\end{cor}
\begin{proof} According to \cite{KaLu} (see
also \cite[Corollary~4.6.12]{vanNbook}), conditions
\eqref{condKaLu1}--\eqref{condKaLu2} imply
$s_0(A)=\omega_0(e^{tA})$. Now
Theorem~\ref{Rebarbthm} gives the result.
\end{proof}

%-------------------------Bibliography-----------------------------


\begin{thebibliography}{99}

\bibitem{Am2} {\sc H. Amann},
{\em Operator-valued Fourier multipliers, vector-valued Besov spaces,
and
applications},  Math.~Nachr.~{\bf 186} (1997), pp. 5--56.


\bibitem{And2} {\sc B.~D.~O.~Anderson}, {\em External and internal
stability of
linear time-varying systems},  SIAM J. Contr. Opt., {\bf 20} no.~3
(1982), pp.~408--413.

\bibitem{BGK} {\sc J.~Ball, I.~Gohberg and M.~A.~Kaashoek}, {\em A
frequency
response function for linear, time-varying systems}, Math.~Control
Signals Systems {\bf 8} (1995), pp.~334--351.

\bibitem{burns} {\sc J.~A.~Burns and B.~B.~King}, {\em A note on the
regularity of
solutions of infinite dimensional Riccati equations},
Appl.~Math.~Lett. Vol.~7, No.~6 (1994), pp.~13--17.

\bibitem{Buse} {\sc C. Buse}, {\em On the
Perron-Bellman theorem for
evolutionary process with exponential
growth in Banach spaces},  New
Zealand J. Math., {\bf 27} (1998), pp. 183--190.

\bibitem{ChiLat} {\sc C. Chicone and Y. Latushkin},
{\em Evolution Semigroups in Dynamical Systems
and Differential Equations},
Mathematical Surveys and Monographs, Vol. 70,
Amer. Math. Soc., Providence, RI, 1999.

\bibitem{Cur2} {\sc R. Curtain}, {\em Equivalence of input-output
stability and
exponential stability for infinite-dimensional systems},
Math. Systems Theory, {\bf 21} (1988) pp.~19--48.

\bibitem{Cur} \sameauthor,
{\em Equivalence of input-output stability and
exponential stability}, Systems and Control Letters, {\bf 12}
(1989), pp.~235--239.

\bibitem{CP} {\sc R.~Curtain and A.~J.~Pritchard},
{\em Infinite Dimensional Linear System Theory},
Lecture Notes in Control and Information Sciences, Vol.~8,
Springer-Verlag, New York, 1978.

\bibitem{CZ} {\sc R.~Curtain and H.~J.~Zwart}, {\em An Introduction to
Infinite-dimensional Linear Systems Theory}, Springer-Verlag,
 New York, 1995.

\bibitem{DK} {\sc J. Daleckij and M. Krein},
{\em Stability of Differential Equations in Banach Space},
Amer. Math. Soc., Providence, RI, 1974.

\bibitem{Datko} {\sc R.~Datko}, {\em Uniform
asymptotic stability of evolutionary
processes in a Banach space},
SIAM J.~Math.~Anal.~{\bf 3} (1972)
pp.~428--445.

\bibitem{FishvN} {\sc A.~Fischer and J.~M.~A.~M.~van Neerven},
{\em Robust stability of $C_0$-semigroups and an application to
stability
of delay equations}, J.~Math.~Anal.~Appl.~{\bf 226} (1998), 82--100.

\bibitem{Hale} {\sc J. Hale}, {\em Ordinary Differential
Equations}, Krieger, 1969.

\bibitem{HIP89} {\sc D.~Hinrichsen, A.~Ilchmann and A.~J.~Pritchard},
{\em Robustness of stability of time-varying linear systems},
J.~Diff.~Eqns., {\bf 82} (1989), pp.~219--250.

\bibitem{HP86} {\sc D.~Hinrichsen and A.~J.~Pritchard}, {\em Stability
radius for structured perturbations and the algebraic Riccati equation},
Systems Control Lett.~{\bf 8} (1986), pp.~105--113.

\bibitem{HP94}
\sameauthor,
{\em Robust stability of
linear evolution operators on Banach spaces},
SIAM J.~Control Optim.,
{\bf 32} no.~6, (1994), pp.~1503--1541.

\bibitem{Huang} {\sc F.~Huang}, {\em Characteristic conditions for
exponential stability
of linear dynamical systems in Hilbert spaces},
Ann.~Diff.~Eqns.~{\bf 1} (1985), pp.~45--53.

\bibitem{JDP} {\sc B.~Jacob, V.~Dragan, A.~J.~Pritchard},
{\em Infinite-dimensional time-varying systems with nonlinear output
feedback.}
Integral Eqns.~Oper.~Theory {\bf 22} (1995),  pp.~440--462.

\bibitem{JacNet} {\sc C.~A.~Jacobson and C.~N.~Nett}, {\em Linear state
space
systems in infinite-dimensional space{\rm:} The role and
characterization of
joint stabilizability/detectability}, IEEE Trans.~Automat.~Control,
{\bf 33}, no.~6 (1988), pp.~541--550.

\bibitem{KaLu} {\sc M.~A.~Kaashoek and S.~M.~Verduyn
Lunel}, {\em An integrability condition for
hyperbolicity of the semigroup}, J.~Diff.~Eqns.~{\bf 112}
(1994), pp.~374--406.

\bibitem{LMS2} {\sc Y.~Latushkin and S.~Montgomery-Smith},
{\em Evolutionary semigroups and Lyapunov theorems in Banach spaces},
J.~Funct.~Anal. {\bf 127} (1995), pp.~173--197.

\bibitem{LMSR} {\sc Y.~Latushkin, S.~Montgomery-Smith,
T.~Randolph}, {\em Evolutionary semigroups and dichotomy of linear
 skew-product flows on locally compact spaces with Banach fibers},
J.~Diff.~Eqns.
 {\bf 125} (1996), pp.~73--116.

\bibitem{LaRa} {\sc Y.~Latushkin and T.~Randolph},
{\em Dichotomy of differential equations on Banach spaces and an
algebra of weighted composition operators},
Integral Equations Operator Theory, {\bf 23} (1995), pp.~472--500.

\bibitem{Log}  {\sc H.~Logemann}, {\em Stabilization and regulation of
infinite-dimensional systems using coprime factorizations,
in Analysis and Optimization of Systems{\rm :} State and Frequency
Domain
Approaches for Infinite-Dimensional Systems} (Sophia-Antipolis, 1992),
edited by R.~F.~Curtain, A.~Bensoussan and J.-L.~Lions,
Lecture Notes in Control and Inform.~Sci., vol.~185, Springer, Berlin
(1993) pp.~102--139.

\bibitem{MiRaSc} {\sc N.~van Minh, R.~R\"abiger and R.~Schnaubelt},
{\em Exponential stability, exponential expansiveness, and exponential
dichotomy of evolution equations on the half-line},
Integral Eqns.~Oper.~Theory, Integral Equations Operator Theory
{\bf 32} (1998), 332--353.

\bibitem{Nagel} {\sc R.~Nagel} (ed.)
{\em One Parameter Semigroups of Positive Operators},
Lecture Notes in Math., no.~1184, Springer-Verlag, Berlin, 1984.

\bibitem{vanN1} {\sc J.~M.~A.~M.~van Neerven},
{\em Characterization of exponential stability of a semigroup of
operators in
terms of its action by convolution on vector-valued function spaces over
$\bbR_+$}, J.~Diff.~Eqns.~{\bf 124} (1996), pp.~324--342.

\bibitem{vanNbook} \sameauthor,
{\em The Asymptotic Behavior of a Semigroup of Linear Operators},
Operator Theory Adv.~Appl. {\bf 88}, Birkhauser, 1996.

\bibitem{Pazy} {\sc A.~Pazy},
{\em Semigroups of Linear Operators and Applications to Partial
Differential Equations},
Springer-Verlag, N.Y./Berlin, 1983.

\bibitem{Phillips} {\sc R. S. Phillips}, {\em Perturbation theory for
semi-groups of linear operators}, Trans. Amer. Math. Soc.
{\bf  74}  (1953), pp. 199--221.

\bibitem{PT89} {\sc A.~J.~Pritchard and S.~Townley}, {\em Robustness of
linear
systems}, J.~Diff. Eqns. {\bf 77} (1989), pp.~254--286.

\bibitem{Prus} {\sc J.~Pr\"{u}ss}, {\em On the spectrum of
$C_0$-semigroups},
Trans.~Amer.~Math.~Soc. {\bf
284} (1984) pp.~847--857.

\bibitem{RS1} {\sc F.~R\"{a}biger and R.~Schnaubelt},
{\em The spectral mapping theorem for evolution semigroups on spaces of
vector-valued functions}, Semigroup Forum {\bf 52} (1996), pp.~225--239.

\bibitem{RRS}  {\sc F. R\"{a}biger, A. Rhandi, R. Schnaubelt},
{\em Perturbation and an
abstract characterization of evolution semigroups},
J.~Math.~Anal.~Appl. {\bf 198} (1996), pp.~516--533.

\bibitem{RRSV} {\sc F.~R\"{a}biger, A.~Rhandi, R.~Schnaubelt and
J.~Voigt}, {\em Non-autonomous Miyadera perturbations},
Differential and Integral Eqns., to appear.

\bibitem{CDC} {\sc T.~Randolph, Y.~Latushkin, S.~Clark},
 {\em Evolution semigroups and stability of time-varying
systems on Banach spaces}, Proceedings of the $36^{th}$
IEEE Conference on Decision and Control, December 1997, pp.~3932--3937.

\bibitem{Rau} {\sc R.~Rau}, {\em Hyperbolic evolution semigroups on
vector valued function spaces}, Semigroup Forum {\bf 48} (1994),
107--118.

\bibitem{Reb1}  {\sc R.~Rebarber}, {\em Conditions for the
 equivalence of internal and external stability for distributed
parameter
systems},
IEEE Trans.~on Automat.~Control vol.~31, no.~6 (1993), 994--998.

\bibitem{Reb2}  \sameauthor,
{\em  Frequency domain methods
for proving the uniform stability of vibrating systems,
in:  Analysis and Optimization of Systems{\rm :} State and
Frequency Domain Approaches for Infinite-Dimensional
Systems} (Sophia-Antipolis, 1992), edited by
R.~F.~Curtain, A.~Bensoussan and J.-L.~Lions, Lecture
Notes in Control and Inform.~Sci., vol.~185, Springer,
Berlin (1993) pp.~366--377.

\bibitem{Renardy} {\sc M.~Renardy}, {\em On the linear
stability of hyperbolic PDEs and viscoelastic flows},
 Z. Angew. Math. Phys. {\bf 45} (1994) pp.~854--865.

\bibitem{RolandDis} {\sc R.~Schnaubelt}, {\em Exponential Bounds and
Hyperbolicity of Evolution Families}, Dissertation,
Eberhard-Karls-Univerist\"at T\"ubingen, 1996.

\bibitem{SK} {\sc R.~Saeks and G.~Knowles}, {\em The Arveson frequency
response and systems theory},  Int.~J.~Control {\bf 42}, no.~3 (1985),
 pp.~639--650.

\bibitem{Weiss1} {\sc G.~Weiss}, {\em Transfer functions of regular
 linear systems, part I{\rm:} Characterizations of regularity},
Trans.~Amer.~Math.~Soc.{\bf 342}, no.~2 (1994), pp.~827--854.

\bibitem{Weiss2} \sameauthor, {\em
Representation of shift invariant
operators
on $L^2$ by $H^\infty$ transfer functions{\rm:} An elementary proof, a
generalization to $L^p$ and a counterexample for $L^\infty$},
 Math.~Control Signals Systems {\bf 4} (1991), pp.~193--203.

\end{thebibliography}
\end{document}